\documentclass[twocolumn]{autart}

\usepackage{algorithm}
\usepackage{algpseudocode}
\usepackage{tikz}
\usetikzlibrary{shapes,arrows}
\usetikzlibrary{arrows.meta}
\tikzset{>={Latex[width=2.5mm,length=2.5mm]}}
\usetikzlibrary{positioning}
\tikzstyle{block}=[draw opacity=0.7,line width=1.4cm]
\tikzset{arrow_e/.style = {->,> = latex'}}
\usepackage{graphicx} 
\usepackage{graphviz}
\usepackage[format=plain, labelsep=period]{caption}
\usepackage{subcaption}
\usepackage{soul}
\usepackage{comment}
\usepackage{amssymb}
\newcommand{\real}{\mathbb{R}}
\usepackage{amsmath}
\newcommand{\vect}[1]{\mathbf{#1}}
\newcommand{\boxend}{\hfill\ensuremath{\Box}}
\newtheorem{assump}{Assumption}
\newcommand{\longthmtitle}[1]{\mbox{}\textit{{(#1):}}}
\usepackage{url}

\begin{document}
\begin{frontmatter}
\title{\large \bf
FORWARD: A Feasible Radial Reconfiguration Algorithm for Multi-Source Distribution Networks\thanksref{footnoteinfo}} 

\thanks[footnoteinfo]{A preliminary version of this paper was presented in 2025 American Control Conference~\cite{JV-RB-SK:24}. }

\author[uci]{Joan Vendrell Gallart}\ead{jvendrel@uci.edu},  
\author[lanl]{Russell Bent}\ead{rbent@lanl.gov}, 
\author[uci]{Solmaz Kia}\ead{solmaz@uci.edu}  

\address[uci]{University of California Irvine, Irvine, California, USA}                             
\address[lanl]{Los Alamos National Laboratory, Los Alamos, New Mexico, USA }   


\begin{abstract}                      
This paper considers an optimal radial reconfiguration problem in multi-source distribution networks, where the goal is to find a radial configuration that minimizes quadratic distribution costs while ensuring all sink demands are met. This problem arises in critical infrastructure systems such as power distribution, water networks, and gas distribution, where radial configurations are essential for operational safety and efficiency. Optimal solution for this problem is known to be NP-hard. In this paper, we prove further that constructing a feasible radial distribution configuration is weakly NP-complete, making exact solution methods computationally intractable for large-scale networks. We propose \texttt{FORWARD} (Feasibility Oriented Random-Walk Inspired Algorithm for Radial Reconfiguration in Distribution Networks), a polynomial-time algorithm that leverages graph-theoretic decomposition and random walk principles to construct feasible radial configurations. Our approach introduces novel techniques including strategic graph partitioning at articulation points, dual graph condensation to address greedy shortsightedness, and capacity-aware edge swapping for infeasibility resolution. We provide rigorous theoretical analysis proving feasibility guarantees and establish a compositional framework enabling parallel processing while preserving optimality properties.
Comprehensive numerical evaluation on networks ranging from IEEE standard test systems to 400-node small-world networks demonstrates that \texttt{FORWARD} consistently outperforms commercial MINLP solvers, achieving optimal or near-optimal solutions in seconds where traditional methods require hours or fail entirely. The algorithm's polynomial-time complexity and scalability make it particularly suitable for real-time distribution network management and as an effective initialization strategy for iterative optimization solvers.
\end{abstract}

\end{frontmatter}

\section{Introduction}
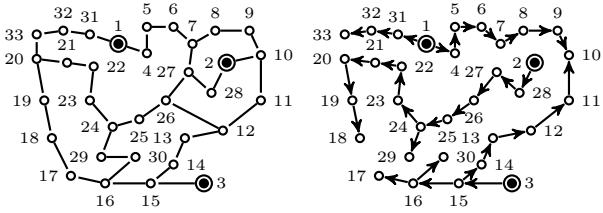
\begin{figure}[t]
    \centering
    \tiny
    \begin{subfigure}{0.2\textwidth}
        \centering
        \begin{tikzpicture}[scale=0.5,>=stealth',shorten >=1pt,auto,node distance=1.5cm,
      thick,main node/.style={circle,fill=black,draw,minimum size=3pt,inner sep=0pt}]
          \node[main node, label=above:$1$] (x1) at (-0.25,0) {};
          \node[main node, fill=none, minimum size=6pt] (x1111) at (-0.25,0) {};
          \node[main node, label=left:$2$] (x2) at (2.6,-0.5) {};
          \node[main node, fill=none, minimum size=6pt] (x2222) at (2.6,-0.5) {};
          \node[main node, label=right:$3$] (x3) at (2,-3.7) {};
          \node[main node, fill=none, minimum size=6pt] (x3333) at (2,-3.7) {};
          \node[main node, label=below:$4$, fill=white] (x4) at (0.5,-0.25) {};
          \node[main node, label=above:$5$, fill=white] (x5) at (0.5,0.45) {};
          \node[main node, label=above:$6$, fill=white] (x6) at (1.2,0.45) {};
          \node[main node, label=above:$7$, fill=white] (x7) at (1.7,0) {};
          \node[main node, label=above:$8$, fill=white] (x8) at (2.3,0.35) {};
          \node[main node, label=above:$9$, fill=white] (x9) at (3.2,0.35) {};
          \node[main node, label=right:$10$, fill=white] (x10) at (3.5,-0.3) {};
          \node[main node, label=right:$11$, fill=white] (x11) at (3.5,-1.5) {};
          \node[main node, label=right:$12$, fill=white] (x12) at (2.5,-2.3) {};
          \node[main node, label=left:$13$, fill=white] (x13) at (1.5,-2.5) {};
          \node[main node, label=right:$14$, fill=white] (x14) at (1.2,-3.2) {};
          \node[main node, label=below:$15$, fill=white] (x15) at (0.6,-3.7) {};
          \node[main node, label=below:$16$, fill=white] (x16) at (-0.6,-3.7) {};
          \node[main node, label=left:$17$, fill=white] (x17) at (-1.5,-3.5) {};
          \node[main node, label=left:$18$, fill=white] (x18) at (-2,-2.5) {};
          \node[main node, label=left:$19$, fill=white] (x19) at (-2.2,-1.5) {};
          \node[main node, label=left:$20$, fill=white] (x20) at (-2.4,-0.4) {};
          \node[main node, label=above:$21$, fill=white] (x21) at (-1.6,-0.5) {};
          \node[main node, label=right:$22$, fill=white] (x22) at (-0.9,-0.6) {};
          \node[main node, label=left:$23$, fill=white] (x23) at (-1,-1.5) {};
          \node[main node, label=left:$24$, fill=white] (x24) at (-0.4,-2.2) {};
          \node[main node, label=below:$25$, fill=white] (x25) at (0.3,-2) {};
          \node[main node, label=below:$26$, fill=white] (x26) at (1,-1.5) {};
          \node[main node, label=left:$27$, fill=white] (x27) at (1.6,-0.75) {};
          \node[main node, label=right:$28$, fill=white] (x28) at (2.2,-1.3) {};
          \node[main node, label=left:$29$, fill=white] (x29) at (-0.7,-3) {};
          \node[main node, label=right:$30$, fill=white] (x30) at (0.2,-3) {};
          \node[main node, label=above:$31$, fill=white] (x31) at (-1,0.25) {};
          \node[main node, label=above:$32$, fill=white] (x32) at (-1.7,0.35) {};
          \node[main node, label=left:$33$, fill=white] (x33) at (-2.4,0.25) {};
          \path[every node/.style={font=\sffamily\small}]
            (x1111) edge (x4)
            (x3333) edge (x15)
            (x4) edge (x5)
            (x5) edge (x6)
            (x6) edge (x7)
            (x7) edge (x8)
            (x8) edge (x9)
            (x9) edge (x10)
            (x10) edge (x11)
            (x11) edge (x12)
            (x12) edge (x13)
            (x13) edge (x14)
            (x14) edge (x15)
            (x15) edge (x16)
            (x16) edge (x30)
            (x16) edge (x17)
            (x2222) edge (x10)
            (x28) edge (x2222)
            (x27) edge (x28)
            (x7) edge (x27)
            (x12) edge (x26)
            (x26) edge (x27)
            (x25) edge (x26)
            (x24) edge (x25)
            (x24) edge (x29)
            (x23) edge (x24)
            (x22) edge (x23)
            (x21) edge (x22)
            (x20) edge (x21)
            (x20) edge (x19)
            (x20) edge (x33)
            (x19) edge (x18)
            (x18) edge (x17)
            (x29) edge (x30)
            (x1111) edge (x31)
            (x31) edge (x32)
            (x32) edge (x33);
    \end{tikzpicture}
        \caption{Original network.}
    \end{subfigure}
    \hspace{3mm}
    \begin{subfigure}{0.2\textwidth}
        \centering
        \begin{tikzpicture}[scale=0.5,>=stealth',shorten <=0.2pt,shorten >=0.2pt,auto,node distance=1.5cm, thick,main node/.style={circle,fill=black,draw,minimum size=3pt,inner sep=0pt}]
          \node[main node, label=above:$1$] (x1) at (-0.25,0) {};
          \node[main node, fill=none, minimum size=6pt] (x1111) at (-0.25,0) {};
          \node[main node, label=left:$2$] (x2) at (2.6,-0.5) {};
          \node[main node, fill=none, minimum size=6pt] (x2222) at (2.6,-0.5) {};
          \node[main node, label=right:$3$] (x3) at (2,-3.7) {};
          \node[main node, fill=none, minimum size=6pt] (x3333) at (2,-3.7) {};
          \node[main node, label=below:$4$, fill=white] (x4) at (0.5,-0.25) {};
          \node[main node, label=above:$5$, fill=white] (x5) at (0.5,0.45) {};
          \node[main node, label=above:$6$, fill=white] (x6) at (1.2,0.45) {};
          \node[main node, label=above:$7$, fill=white] (x7) at (1.7,0) {};
          \node[main node, label=above:$8$, fill=white] (x8) at (2.3,0.35) {};
          \node[main node, label=above:$9$, fill=white] (x9) at (3.2,0.35) {};
          \node[main node, label=right:$10$, fill=white] (x10) at (3.5,-0.3) {};
          \node[main node, label=right:$11$, fill=white] (x11) at (3.5,-1.5) {};
          \node[main node, label=right:$12$, fill=white] (x12) at (2.5,-2.3) {};
          \node[main node, label=left:$13$, fill=white] (x13) at (1.5,-2.5) {};
          \node[main node, label=right:$14$, fill=white] (x14) at (1.2,-3.2) {};
          \node[main node, label=below:$15$, fill=white] (x15) at (0.6,-3.7) {};
          \node[main node, label=below:$16$, fill=white] (x16) at (-0.6,-3.7) {};
          \node[main node, label=left:$17$, fill=white] (x17) at (-1.5,-3.5) {};
          \node[main node, label=left:$18$, fill=white] (x18) at (-2,-2.5) {};
          \node[main node, label=left:$19$, fill=white] (x19) at (-2.2,-1.5) {};
          \node[main node, label=left:$20$, fill=white] (x20) at (-2.4,-0.4) {};
          \node[main node, label=above:$21$, fill=white] (x21) at (-1.6,-0.5) {};
          \node[main node, label=right:$22$, fill=white] (x22) at (-0.9,-0.6) {};
          \node[main node, label=left:$23$, fill=white] (x23) at (-1,-1.5) {};
          \node[main node, label=left:$24$, fill=white] (x24) at (-0.4,-2.2) {};
          \node[main node, label=below:$25$, fill=white] (x25) at (0.3,-2) {};
          \node[main node, label=below:$26$, fill=white] (x26) at (1,-1.5) {};
          \node[main node, label=left:$27$, fill=white] (x27) at (1.6,-0.75) {};
          \node[main node, label=right:$28$, fill=white] (x28) at (2.2,-1.3) {};
          \node[main node, label=left:$29$, fill=white] (x29) at (-0.7,-3) {};
          \node[main node, label=right:$30$, fill=white] (x30) at (0.2,-3) {};
          \node[main node, label=above:$31$, fill=white] (x31) at (-1,0.25) {};
          \node[main node, label=above:$32$, fill=white] (x32) at (-1.7,0.35) {};
          \node[main node, label=left:$33$, fill=white] (x33) at (-2.4,0.25) {};
          \path[every node/.style={font=\sffamily\small}]
            (x1111) edge[->] (x4)
            (x3333) edge[->] (x15)
            (x4) edge[->] (x5)
            (x5) edge[->] (x6)
            (x6) edge[->] (x7)
            (x7) edge[->] (x8)
            (x8) edge[->] (x9)
            (x9) edge[->] (x10)
            (x11) edge[->] (x10)
            (x12) edge[->] (x11)
            (x13) edge[->] (x12)
            (x14) edge[->] (x13)
            (x15) edge[->] (x14)
            (x15) edge[->] (x16)
            (x16) edge[->] (x30)
            (x16) edge[->] (x17)
            (x2222) edge[->] (x28)
            (x28) edge[->] (x27)
            (x27) edge[->] (x26)
            (x26) edge[->] (x25)
            (x25) edge[->] (x24)
            (x24) edge[->] (x29)
            (x24) edge[->] (x23)
            (x23) edge[->] (x22)
            (x22) edge[->] (x21)
            (x21) edge[->] (x20)
            (x20) edge[->] (x19)
            (x19) edge[->] (x18)
            (x1111) edge[->] (x31)
            (x31) edge[->] (x32)
            (x32) edge[->] (x33);
    \end{tikzpicture}
        \caption{Radial configuration.}
    \end{subfigure}
    \caption{{\small  In the optimal reconfiguration problem, the highlighted nodes in dark are the sources and the remaining nodes are the sinks. In radial configuration, some nodes, e.g., sink node 10 may receive receive input from two different edge; despite that there is no cycle in the graph. The network used here is the IEEE 33 network~\cite{ieee33}. }}  
    \label{fig::network_example} 
\end{figure}

\vspace{-0.1in}
In today's rapidly evolving technological and infrastructural landscape, distribution networks, encompassing electricity, gas, and water irrigation systems, are fundamental to resource delivery across diverse geographical areas. In these distribution networks, there is often a need for radial reconfiguration (see Fig.~\ref{fig::network_example}), driven by the complex physical properties and behaviors inherent in resource flow systems, as well as the necessity to adhere to engineering and safety standards. Radial configurations help minimize losses and ensure system stability, which are essential for both economic viability and environmental sustainability~\cite{AM-HB:75,MEB-FFW:89,RAJ-RS-BCP:12}.
With the increasing complexity of modern distribution systems, there is a need for computationally fast solutions  to obtain  reliable resource distribution.

Optimizing distribution over all the possible radial reconfiguration often is cast as  Mixed-Integer Non-Linear Programming (MINLP) problems~\cite{RAJ-RS-BCP:12,AFZ-AV:14,RYY-TP-RR-LHM:24}. The complexity of solving these optimization problems arises from the exponential number of possible configurations, particularly as the network scales, making it an NP-hard problem~\cite{NS:19}. Traditional MINLP methods often fall short of delivering an optimal solution due to their computational demands, especially in large-scale networks. Recent advancements have shifted towards suboptimal algorithms that can efficiently handle this complexity. For instance, leveraging graph theory and random walk processes has shown promise in simplifying the search space and reducing computational effort~\cite{EP-CB-JV:20}.

\setlength{\textfloatsep}{3pt}

This paper constructs a polynomial-time  feasible solution for a class of optimal radial reconfiguration problem in a distribution network with multiple sources and sinks. Our method, termed as \texttt{FORWARD}: \textbf{F}easibility \textbf{O}riented \textbf{R}andom-\textbf{Wa}lk Inspired Algorithm for \textbf{R}adial Reconfiguration in \textbf{D}istribution Networks, leverages the graph-like characteristics of distribution systems and takes into account how potential flows naturally navigate through a network. Our algorithm uses a greedy radial construction process starting from source nodes, incrementally adding edges while considering the flow across constructed components and the demanded output of the remaining nodes to reach. We draw inspiration from the similarities between electric flow in power networks and random walks~\cite{randomwalk}, developing a novel `sampling' method\footnote{Although our radial configuration construction is deterministic, we use `sampling' as a conceptual analogy.} for constructing radial configurations. In the random walk approach to describe electricity distribution in a given network, a weight proportional to the inverse of the resistance along the corresponding link in the power network is assigned to each link. This creates a notion of the edge weights being the conductance of the edge. Just as electricity flows through paths of least resistance, a random walk probabilistically selects paths based on transition probabilities influenced by edge weights.  \texttt{FORWARD} respects capacity constraints and ensures the network remains feasible and efficient by reducing the sample space to feasible edges and prioritizing critical edges. 

We introduce an innovative mechanism (\textsf{Net-Concad} function explained in Section~\ref{sec:net-concad}) in our incremental `sampling' process that addresses the shortsightedness of greedy selections by informing the process of the comprehensive demand of the remaining nodes. We carry out a rigorous evaluation of the algorithm, including its formal feasibility guarantees.
Extensive numerical experiments demonstrate the algorithm's efficacy across various distribution networks, highlighting its potential for real-world applications~\cite{russell}. 

In summary, the main contributions of this work:
\begin{itemize}
    \item We prove that the problem of constructing a feasible radial configuration for multi-source distribution networks is weakly NP-complete, a new theoretical result that highlights the computational intractability of the problem and motivates the need for efficient algorithms.
    \item We propose a novel polynomial-time algorithm, named \texttt{FORWARD}, that leverages a compositional framework to guarantee feasible solutions and rapidly find high-quality solutions with a time complexity of $\mathcal{O}(n^2 \log n)$ on sparse networks.
    \item We introduce several key algorithmic innovations, including a dual graph condensation technique to overcome greedy limitations and a capacity-aware edge swapping method to resolve infeasibilities.
    \item We demonstrate that \texttt{FORWARD} consistently outperforms commercial MINLP solvers on large-scale networks (up to 400 nodes), achieving optimal or near-optimal solutions in seconds where other methods fail or take hours, making it highly suitable for real-time network management.
\end{itemize}

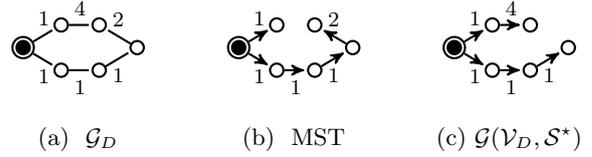
\begin{figure}[t]
    \centering
        \begin{subfigure}{0.12\textwidth}
            \centering
            \begin{tikzpicture}[>=stealth',shorten      >=1pt,auto,node distance=1.5cm,
            thick,main node/.style={circle,fill=black,draw,minimum size=5pt,inner sep=0pt}]
              \node[main node, fill=black] (x1) at (0,0) {};
              \node[main node, fill=none, minimum size=8pt] (x1111) at (0,0) {};
              \node[main node, fill=white] (x2) at (0.5,-0.3) {};
              \node[main node, fill=white] (x3) at (0.5,0.3) {};
              \node[main node, fill=white] (x4) at (1,-0.3) {};
              \node[main node, fill=white] (x5) at (1.5,0) {};
              \node[main node, fill=white] (x6) at (1,0.3) {};
                \draw  (x1111) -- node[below] {\scriptsize 1} (x2);
                \draw  (x1111) -- node[above] {\scriptsize 1} (x3);
                \draw  (x2) -- node[below] {\scriptsize 1} (x4);
                \draw  (x4) -- node[below] {\scriptsize 1} (x5);
                \draw (x5) -- node[above] {\scriptsize 2} (x6);
                \draw (x3) -- node[above] {\scriptsize 4} (x6);
            \end{tikzpicture}
            \caption{{ $\mathcal{G}_D$}}
        \end{subfigure}
        \hspace{5mm}
        \begin{subfigure}{0.12\textwidth}
            \centering
            \begin{tikzpicture}[>=stealth',shorten      >=1pt,auto,node distance=1.5cm,
            thick,main node/.style={circle,fill=black,draw,minimum size=5pt,inner sep=0pt}]
              \node[main node, fill=black] (x1) at (0,0) {};
              \node[main node, fill=none, minimum size=8pt] (x1111) at (0,0) {};
              \node[main node, fill=white] (x2) at (0.5,-0.3) {};
              \node[main node, fill=white] (x3) at (0.5,0.3) {};
              \node[main node, fill=white] (x4) at (1,-0.3) {};
              \node[main node, fill=white] (x5) at (1.5,0) {};
              \node[main node, fill=white] (x6) at (1,0.3) {};
                \draw[->] (x1111) -- node[below] {\scriptsize 1} (x2);
                \draw[->] (x1111) -- node[above] {\scriptsize 1} (x3);
                \draw[->] (x2) -- node[below] {\scriptsize 1} (x4);
                \draw[->] (x4) -- node[below] {\scriptsize 1} (x5);
                \draw[->] (x5) -- node[above] {\scriptsize 2} (x6);
            \end{tikzpicture}
            \caption{{ MST}}
        \end{subfigure}
        \hspace{5mm}
        \begin{subfigure}{0.12\textwidth}
            \centering
            \begin{tikzpicture}[>=stealth',shorten      >=1pt,auto,node distance=1.5cm,
            thick,main node/.style={circle,fill=black,draw,minimum size=5pt,inner sep=0pt}]
              \node[main node, fill=black] (x1) at (0,0) {};
              \node[main node, fill=none, minimum size=8pt] (x1111) at (0,0) {};
              \node[main node, fill=white] (x2) at (0.5,-0.3) {};
              \node[main node, fill=white] (x3) at (0.5,0.3) {};
              \node[main node, fill=white] (x4) at (1,-0.3) {};
              \node[main node, fill=white] (x5) at (1.5,0) {};
              \node[main node, fill=white] (x6) at (1,0.3) {};
                \draw[->] (x1111) -- node[below] {\scriptsize 1} (x2);
                \draw[->] (x1111) -- node[above] {\scriptsize 1} (x3);
                \draw[->] (x2) -- node[below] {\scriptsize 1} (x4);
                \draw[->] (x4) -- node[below] {\scriptsize 1} (x5);
                \draw[->] (x3) -- node[above] {\scriptsize 4} (x6);
            \end{tikzpicture}
            \caption{{$\mathcal{G}(\mathcal{V}_D,\mathcal{S}^\star)$}}
        \end{subfigure}
    \caption{{\footnotesize Example where radial distribution constructed from MST (plot (b)) is not the minimum radial configuration, $\mathcal{G}(\mathcal{V}_D,\mathcal{S}^\star)$, (plot(c)). This phenomenon is due to the quadratic nature of the cost; let the demand at each sink be $d$. The source node, highlighted in bold, can supply input $5d$. The edge weights are shown on the edges. In the feasible radial distribution network (b) the cost is $1\cdot(d)^2 +1\cdot(4d)^2+1\cdot(3d)^2+1\cdot(2d)^2+2\cdot(d)^2 =32d^2$, meanwhile in the feasible radial distribution network (c) the cost is $1\cdot(2d)^2 +4\cdot(1d)^2  +1\cdot(3d)^2+1\cdot(2d)^2+1\cdot(1d)^2=22d^2$.}}
\label{fig::quadratic_effect}
\end{figure}

\begin{figure}[t]
    \centering
        \begin{subfigure}{0.1\textwidth}
           \centering
            \begin{tikzpicture}[>=stealth',shorten      >=1pt,auto,node distance=1.5cm,
            thick,main node/.style={circle,fill=black,draw,minimum size=5pt,inner sep=0pt}]
              \node[main node, fill=black] (x1) at (0,0) {};
              \node[main node, fill=none, minimum size=8pt] (x1111) at (0,0) {};
              \node[main node, fill=none, minimum size=8pt] (x6666) at (1.5,0) {};
              \node[main node, fill=white] (x3) at (0.5,0) {};
              \node[main node, fill=white] (x4) at (1,0) {};
              \node[main node, fill=black] (x6) at (1.5,0) {};
              \node[main node, fill=white] (x7) at (0.5,-0.7) {};
              \node[main node, fill=white] (x8) at (1,-0.7) {};
                \draw (x1111) -- node[above] {\scriptsize 1} (x3);
                \draw (x3) -- node[above] {\scriptsize 1} (x4);
                \draw (x3) -- node[left] {\scriptsize 1} (x7);
                \draw (x4) -- node[above] {\scriptsize 1} (x6666);
                \draw (x7) -- node[above] {\scriptsize 1} (x8);
                \draw (x4) -- node[right] {\scriptsize 1.5} (x8);
            \end{tikzpicture}
            \caption{{ $\mathcal{G}$}}
        \end{subfigure}
    \hspace{5mm}
        \begin{subfigure}{0.1\textwidth}
        \centering
            \begin{tikzpicture}[>=stealth',shorten      >=1pt,auto,node distance=1.5cm,
            thick,main node/.style={circle,fill=black,draw,minimum size=5pt,inner sep=0pt}]
              \node[main node, fill=black] (x1) at (0,0) {};
              \node[main node, fill=none, minimum size=8pt] (x1111) at (0,0) {};
              \node[main node, fill=none, minimum size=8pt] (x6666) at (1.7,0) {};
              \node[main node, fill=white] (x3) at (0.6,0) {};
              \node[main node, fill=white] (x5) at (1.1,0) {};
              \node[main node, fill=black] (x6) at (1.7,0) {};
              \node[main node, fill=white] (x7) at (0.6,-0.7) {};
              \node[main node, fill=white] (x8) at (1.1,-0.7) {};
                \draw (x1111) -- node[above] {\scriptsize 1} (x3);
                \draw (x6666) -- node[above] {\scriptsize 1} (x5);
                \draw (x3) -- node[left] {\scriptsize 1} (x7);
                \draw (x5) -- node[above] {\scriptsize 1} (x3);
                \draw (x7) -- node[above] {\scriptsize 1} (x8);
            \end{tikzpicture}
            \caption{{ MST}}
        \end{subfigure}
        \hspace{5mm}
        \begin{subfigure}{0.1\textwidth}
        \centering
            \begin{tikzpicture}[>=stealth',shorten      >=1pt,auto,node distance=1.5cm,
            thick,main node/.style={circle,fill=black,draw,minimum size=5pt,inner sep=0pt}]
              \node[main node, fill=black] (x1) at (0,0) {};
              \node[main node, fill=none, minimum size=8pt] (x1111) at (0,0) {};
              \node[main node, fill=none, minimum size=8pt] (x6666) at (1.7,0) {};
              \node[main node, fill=white] (x3) at (0.6,0) {};
              \node[main node, fill=white] (x4) at (1.1,0) {};
              \node[main node, fill=black] (x6) at (1.7,0) {};
              \node[main node, fill=white] (x7) at (0.6,-0.7) {};
              \node[main node, fill=white] (x8) at (1.1,-0.7) {};
                \draw (x1111) -- node[above] {\scriptsize 1} (x3);
                \draw (x3) -- node[left] {\scriptsize 1} (x7);
                \draw (x6666) -- node[above] {\scriptsize 1} (x4);
                \draw (x4) -- node[right] {\scriptsize 1.5} (x8);
            \end{tikzpicture}
            \caption{{ MSF}}
        \end{subfigure}
    \caption{{\small An example where MSF results is a better outcome than MST: sources highlighted provide an input of $2d$ each and the sink nodes demand each $d$. }}\vspace{0.1in}
    \label{fig::complexity_msf}
\end{figure}
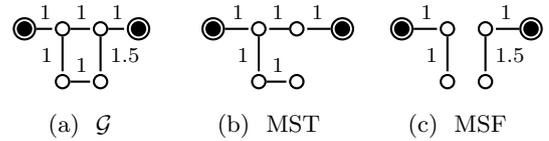

\emph{Related work}: Several methods in the literature, especially in power networks applications, leverage topological properties to address the reconfiguration problem \cite{clark}. For example, \cite{AK-GY-SB-EN-MCC-TL-EP:18} \cite{spectral} use spectral clustering followed by local greedy search to identify radial configurations. In \cite{cuts}, the reconfiguration problem is linked to the maximum flow problem, based on Ford-Fulkerson’s work \cite{frodfulkerson}, which has been extensively studied in both single- and multi-source contexts \cite{maxflow}, \cite{maxflow_multi}, \cite{maxflow_multi2}. These approaches often include a repair procedure for unfeasible solutions, which eventually results inefficient in practice. While some methods address minimum-cost distribution, radiality remains a critical constraint \cite{non_radial}, and cycle-breaking methods used to tackle it offer limited guarantees for large graphs \cite{LC-RK-YPL-RP-MPG-SS:22}. Other approaches, such as those based on the minimum spanning tree (MST) problem, face challenges due to the constraints of the reconfiguration problem, making the MST approach NP-hard in this context (see Fig.~\ref{fig::quadratic_effect}). Moreover, although algorithms like Kruskal's or Prim's \cite{tate2016prim} can find MSTs in polynomial time, they lose their optimality when multiple sources are involved, note that the minimum spanning forest (MSF) is not necessarily a subset of the MST (see Fig.~\ref{fig::complexity_msf}).

\section{Problem Setting}
\label{sec::intro}
\vspace{-0.1in}
We consider a network-flow distribution problem over a bidirectional distribution network $\mathcal{G}_D=\mathcal{G}(\mathcal{V}_D, \mathcal{E}_D)$ with $ |\mathcal{V}_D| = N $ nodes, $|\mathcal{E}_D| = m$ edges, a set of $n_g$ source nodes $\mathcal{V}_g \subset \mathcal{V}_D$, and $n_c = N - n_g $ sink nodes $\mathcal{V}_c = \mathcal{V}_D \backslash \mathcal{V}_g$, each with specified input and output.  We let $\vect{d}\in\real_{\geq0}^N$ be the output vector whose non-zero elements correspond to sink nodes and $\vect{g}\in\real_{\geq 0}^N$ be the input vector whose non-zero elements belong to source nodes; i.e., $g_i\in\real_{>0}$ (resp. $d_j\in\real_{>0}$) for $i\in\mathcal{V}_g$ (resp. $j\in\mathcal{V}_c$) and $ g_i=0$ (resp. $d_j=0$) otherwise.  The assumption is that inputs match outputs, i.e., $\sum_{i\in\mathcal{V}_g} g_i=\sum_{i\in\mathcal{V}_c} d_i$. 

The goal is to find an (oriented) \emph{radial configuration} that delivers the input flow from the source nodes to the sink nodes with minimal overall cost to meet the output demands. An example valid radial configuration is shown in Fig.~\ref{fig::network_example}, which shows a constructed radial configuration (plot (b)) from a given distribution network (plot (a)). This cost results from the `resistance' or `toll' along the edges, characterized as a quadratic function of the flow across them.

\begin{defn}
   \textbf{(Set of radial configurations)}
   \label{def::radial}
    A radial configuration is a polyforest\footnote{A polyforest (or directed forest or oriented forest) is a directed acyclic graph whose underlying undirected graph is a forest. A forest is a type of graph that contains no loops. Consequently, forests consist solely of trees that might be disconnected, leading to the term `forest' being used \cite{forests}.} that includes all the nodes $\mathcal{V}_D$, has roots at $\mathcal{V}_g$ and the undirected version of its edges are subset of $\mathcal{E}_D$. We denote the set of these polyforest digraphs by $\mathcal{F}(\mathcal{G}_D,\mathcal{V}_g)$. For brevity, when clear from context, we will use only $\mathcal{F}$. \boxend
\end{defn}

The problem of interest is formalized as 
\vspace{-0.08in}
\begin{subequations}\label{eqn::problem1}
\begin{align}
\min ~&\sum\nolimits_{(i,j)\in\mathcal{S}} C_{i,j}\cdot x_{i,j}^2,\quad \text{subject to}\label{eqn::distribution_cost}\\
     &~~ \mathcal{G}(\mathcal{V}_D,\mathcal{S})\in\mathcal{F},\label{eqn::distribution_cost_radial}\\
     &~~0\leq\vect{x}(\mathcal{S})\leq\bar{\vect{x}}(\mathcal{S}) ,\label{eqn::distribution_cost_capacity}\\
     &~~A(\mathcal{S})\, \vect{x}(\mathcal{S}) ={\vect{g}} - {\vect{d}}, \label{eqn::problem1-kirchof}
\end{align}
\end{subequations}
where $x_{i,j}\in\real_{>0}$ is the flow across the link $(i,j)$, $A(\mathcal{S})$ is the incidence matrix of the radial configuration $\mathcal{S}$ (a decision variable of the optimization problem), $C_{i,j}$ is the coefficient of the cost of edge $(i,j)\in\mathcal{S}$. The constraint~\eqref{eqn::distribution_cost_radial} confines the radial configuration to $\mathcal{F}$, \eqref{eqn::distribution_cost_capacity} enforce capacity across each edge in $\mathcal{S}$, and \eqref{eqn::problem1-kirchof} enforces the flow conservation (Kirchhoff's law) at the nodes. We consider the following assumption across this paper.
\begin{assump}\longthmtitle{Feasibility and feasible radial distribution configuration}\label{assump::feas}
    The feasible solution set of the optimization problem~\eqref{eqn::problem1} is non-empty. We refer to any feasible solution as a \emph{feasible radial distribution configuration}.  \boxend
\end{assump}
This assumption means that at least there is one feasible radial distribution configuration, for which the inputs can meet the specified output, despite the capacity~bounds.

Problem~\eqref{eqn::problem1} is highly relevant to various potential flow applications, including natural gas, water, and electricity distribution networks. For instance, in power systems, the efficient and reliable distribution of energy has become increasingly critical with the growing integration of renewable energy and distributed generation sources. Modern power distribution networks, comprising multiple distributed generators, must operate in a radial configuration to adhere to engineering and safety standards. Minimizing energy loss within these networks is vital for both economic viability and environmental sustainability. Consequently, radial configurations cannot be chosen arbitrarily; they must be designed to ensure feasibility and to optimize energy loss (cost of operation).

\section{Notations}
\vspace{-0.1in}
In this section we introduce our essential notation and review the definition of some graph theoretic and algorithmic concepts and tools we use in the paper: 
\begin{itemize}
    \item $\mathcal{S}$: the set of sampled directed edges in the FORWARD algorithm.
    \item $\mathcal{V}(\mathcal{S})$: the nodes of the sampled edge set $\mathcal{S}$.
    \item $\mathcal{E}(\mathcal{V}_i)$: the set of edges in $\mathcal{E}_D$ interconnecting nodes in $\mathcal{V}_i$.
    \item $\mathcal{N}(\mathcal{V}_i)$: the set of nodes in $\mathcal{V}_D\backslash\mathcal{V}_i$ directly connected to nodes in $\mathcal{V}_i$.
    \item $(i\to j)$: a directed edge indicating flow from node $i$ to $j$.
\end{itemize}
We define \emph{current nodal value} vector (often referred to as value vector) $\vect{p}$, associated with the nodes of $\mathcal{V}_D$ as the unified input/output value vector across the distribution network, initialize at $\vect{p}_i=-d_i<0$ for $i\in\mathcal{V}_c$ and $\vect{p}_i=g_i>0$ for $i\in\mathcal{V}_g$. In our algorithmic solutions, this vector can be updated at each iteration of the algorithm, depending on the process described. 

A \emph{quasi-bipartite graph} is a graph whose vertices can be partitioned into two disjoint subsets, say $U$ and $V$, such that the vast majority of edges connect vertices from $U$ to $V$. Unlike strictly bipartite graphs, however, a small number of edges may exist within $U$ or within $V$~\cite{DH:22}. A \emph{$2$-core} subgraph, or $2$-core, is a maximal subgraph of a graph where all vertices have a degree of at least $2$. You can find the $2$-core of a graph $G$ by repeatedly removing all vertices that have a degree less than $2$ in $G$. When you remove a vertex, its incident edges are also removed, which can reduce the degrees of its neighbors. This process is iterated until no more vertices can be removed. The remaining subgraph is the $2$-core. An \emph{articulation} node in a graph is a node which, if removed, causes the graph to be disconnected. A bi-connected graph is a graph with no articulation nodes. In a graph, a \emph{pendent node} (also called a leaf node or end vertex) is a node (or vertex) that has a degree of exactly one, meaning it is connected to only one other node. A \emph{pendent edge} is an edge connected to a pendent node. We let $\mathcal{G}_P$ to denote the 2-core subgraph of $\mathcal{G}_D$ and $\bar{L}$ is the total number of articulation nodes in $\mathcal{G}_P$.

A \emph{queue} is a linear data structure in which elements are \emph{inserted} at the rear and \emph{retrieved} from the front (FIFO data structure: First-In, First-Out).

\section{Hardness analysis and objective statement}
\label{sec::complexity}
\vspace{-0.1in}
Problem~\eqref{eqn::problem1} is known to be NP-hard to solve optimally~\cite{AM-HB:75}. Since this paper aims to construct a feasible suboptimal solution, we next investigate the computational complexity of finding any feasible solution for Problem~\eqref{eqn::problem1}. Our analysis, formalized in Theorem~\ref{thm::hardness}, reveals that determining a feasible distribution configuration for problem~\eqref{eqn::problem1} is weakly NP-complete. This classification is a particularly significant finding, as it indicates that despite its overall NP-hardness, the problem's dependence on the numerical values of its input permits the development of pseudo-polynomial time~algorithms to find a feasible solution.

\begin{figure}[t]
    \centering
    \begin{subfigure}{0.10\textwidth}
        \centering
    \begin{tikzpicture}[>=stealth',shorten      >=1pt,auto,node distance=1.5cm,
           thick,main node/.style={circle,fill=black,draw,minimum size=8pt,inner sep=0pt}]
         \node[main node, fill=black, label=left:{\scriptsize $9$}] (x1) at (0,0) {};
         \node[main node, fill=none, minimum size=11pt] (x1111) at (0,0) {{{\tiny{\color{white}$s_1$}}}};
         \node[main node, fill=black, label=right:{\scriptsize $9$}] (x2) at (1.5,0) {};
         \node[main node, fill=none, minimum size=11pt] (x2222) at (1.5,0) {{{\tiny{\color{white}$s_2$}}}};
         \node[main node, fill=white, label=above:{\scriptsize $4$}] (x3) at (0.75,1.5) {{\tiny$v_4$}};
         \node[main node, fill=white, label=above:{\scriptsize $5$}] (x4) at (0.75,0.75) {{\tiny$v_3$}};
         \node[main node, fill=white, label=above:{\scriptsize $4$}] (x5) at (0.75,0) {{\tiny$v_2$}};
         \node[main node, fill=white, label=above:{\scriptsize $3$}] (x6) at (0.75,-0.75) {{\tiny$v_1$}};
         \node[main node, fill=white, label=above:{\scriptsize $2$}] (x7) at (0.75,-1.5) {{\tiny $v_0$}};
         \draw  (x1111) -- node[above] {\tiny \textit{4}} (x3);
         \draw  (x2222) -- node[above] {\tiny \textit{4}} (x3);
         \draw  (x1111) -- node[above] {\tiny \textit{5}} (x4);
         \draw  (x2222) -- node[above] {\tiny \textit{5}} (x4);
         \draw  (x1111) -- node[above] {\tiny \textit{4}} (x5);
         \draw  (x2222) -- node[above] {\tiny \textit{4}} (x5);
         \draw  (x1111) -- node[below] {\tiny \textit{3}} (x6);
         \draw  (x2222) -- node[below] {\tiny \textit{3}} (x6);
         \draw  (x1111) -- node[below] {\tiny \textit{1}} (x7);
         \draw  (x2222) -- node[below] {\tiny \textit{1}} (x7);
        \end{tikzpicture}
        \caption{{Distribution network $\mathcal{G}_D$.}}
   \end{subfigure}
   \hspace{6mm}
    \begin{subfigure}{0.10\textwidth}
        \centering
        \begin{tikzpicture}[>=stealth',shorten      >=1pt,auto,node distance=1.5cm,
           thick,main node/.style={circle,fill=black,draw,minimum size=8pt,inner sep=0pt}]
         \node[main node, fill=black, label=left:{\scriptsize $0$}] (x1) at (0,0) {};
         \node[main node, fill=none, minimum size=11pt] (x1111) at (0,0) {{{\tiny{\color{white}$s_1$}}}};
         \node[main node, fill=black, label=right:{\scriptsize $0$}] (x2) at (1.5,0) {};
         \node[main node, fill=none, minimum size=11pt] (x2222) at (1.5,0) {{{\tiny{\color{white}$s_2$}}}};
         \node[main node, fill=white, label=above:{\scriptsize $0$}] (x3) at (0.75,1.5) {{\tiny$v_4$}};
         \node[main node, fill=white, label=above:{\scriptsize $0$}] (x4) at (0.75,0.75) {{\tiny$v_3$}};
         \node[main node, fill=white, label=above:{\scriptsize $0$}] (x5) at (0.75,0) {{\tiny$v_2$}};
         \node[main node, fill=white, label=above:{\scriptsize $0$}] (x6) at (0.75,-0.75) {{\tiny$v_1$}};
         \node[main node, fill=white, label=above:{\scriptsize $0$}] (x7) at (0.75,-1.5) {{\tiny$v_0$}};
         \draw[->]  (x1111) -- node[above] {\tiny \textit{4}} (x3);
         \draw[dashed, color=gray]  (x2222) -- node[above] {\tiny \textit{4}} (x3);
         \draw[dashed, color=gray]   (x1111) -- node[above] {\tiny \textit{5}} (x4);
         \draw[->]  (x2222) -- node[above] {\tiny \textit{5}} (x4);
         \draw[->]  (x1111) -- node[above] {\tiny \textit{4}} (x5);
         \draw [dashed, color=gray]  (x2222) -- node[above] {\tiny \textit{4}} (x5);
         \draw[dashed, color=gray]   (x1111) -- node[below] {\tiny \textit{3}} (x6);
         \draw[->]  (x2222) -- node[below] {\tiny \textit{3}} (x6);
         \draw[->]  (x1111) -- node[below] {\tiny \textit{1}} (x7);
         \draw[->]  (x2222) -- node[below] {\tiny \textit{1}} (x7);
        \end{tikzpicture}
        \caption{{A feasible radial distribution.}}
   \end{subfigure}
   \hspace{6mm}
    \begin{subfigure}{0.10\textwidth}
        \centering
        \begin{tikzpicture}[>=stealth',shorten      >=1pt,auto,node distance=1.5cm,
           thick,main node/.style={circle,fill=black,draw,minimum size=7pt,inner sep=0pt}]
         \node[main node, fill=black, label=left:{\scriptsize $0$}] (x1) at (0,0) {{{\tiny{\color{white}$s_1$}}}};
         \node[main node, fill=none, minimum size=11pt] (x1111) at (0,0) {};
         \node[main node, fill=black, label=right:{\scriptsize $0$}] (x2) at (1.5,0) {{{\tiny{\color{white}$s_2$}}}};
         \node[main node, fill=none, minimum size=11pt] (x2222) at (1.5,0) {};
         \node[main node, fill=white, label=above:{\scriptsize $0$}] (x3) at (0.75,1.5) {{\tiny$v_4$}};
         \node[main node, fill=white, label=above:{\scriptsize $0$}] (x4) at (0.75,0.75) {{\tiny$v_3$}};
         \node[main node, fill=white, label=above:{\scriptsize $0$}] (x5) at (0.75,0) {{\tiny$v_2$}};
         \node[main node, fill=white, label=above:{\scriptsize $0$}] (x6) at (0.75,-0.75) {{\tiny$v_1$}};
         \node[main node, fill=white, label=above:{\scriptsize $0$}] (x7) at (0.75,-1.5) {{\tiny$v_0$}};
         \draw[->]  (x1111) -- node[above] {\tiny \textit{4}} (x3);
         \draw[->]  (x2222) -- node[above] {\tiny \textit{4}} (x3);
         \draw[dashed, color=gray]   (x1111) -- node[above] {\tiny \textit{5}} (x4);
         \draw[->]  (x2222) -- node[above] {\tiny \textit{5}} (x4);
         \draw[->]  (x1111) -- node[above] {\tiny \textit{4}} (x5);
         \draw [dashed, color=gray]  (x2222) -- node[above] {\tiny \textit{4}} (x5);
         \draw[->]   (x1111) -- node[below] {\tiny \textit{3}} (x6);
         \draw[->]  (x2222) -- node[below] {\tiny \textit{3}} (x6);
         \draw[->]  (x1111) -- node[below] {\tiny \textit{1}} (x7);
         \draw[->]  (x2222) -- node[below] {\tiny \textit{1}} (x7);
        \end{tikzpicture}
        \caption{{A non-radial distribution.}}
   \end{subfigure}
    \caption{{\footnotesize An example of a distribution scenario described in the proof of Theorem~\ref{thm::hardness}. Capacities are shown as edge weights. Source nodes are shown in dark }}\vspace{0.1in}
    \label{fig::np_complete}
\end{figure}
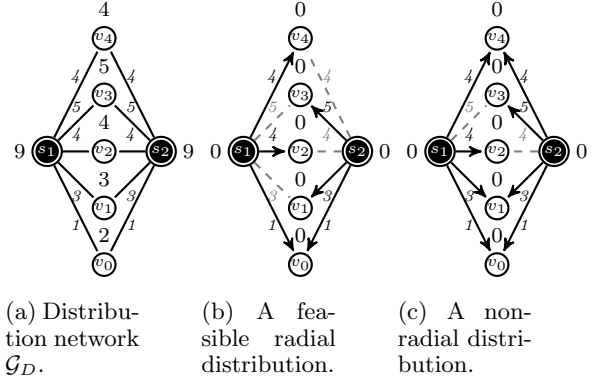
\begin{thm}
    \longthmtitle{Computational complexity of finding a feasible radial distribution
configuration for problem~\eqref{eqn::problem1}}
    Finding a feasible radial distribution  configuration for problem~\eqref{eqn::problem1} is weakly NP-complete.
    \label{thm::hardness}
\end{thm}
\begin{pf}
To prove the problem is in NP, we must demonstrate that a given candidate solution can be verified in polynomial time. A candidate solution consists of a set of chosen directed edges $\mathcal{S} \subseteq \mathcal{E}_D$ (directed graph $(\mathcal{V}_D,\mathcal{S})$) and their corresponding flow values $\vect{x}(\mathcal{S})$. Recall $N = |\mathcal{V}_D|$ and $m = |\mathcal{E}_D|$. We verify the solution in three steps with polynomial computational cost: 
\begin{itemize}
    \item \emph{Step 1 (Check capacity constraint is respected; computational complexity $O(m)$):}  Check whether the capacity constraint along very edge in $\mathcal{S}$ is satisfied. 

    \item \emph{Step 2 (Check flow conservation at each node; computational complexity $O(m+n)$)}: Check whether flow conservation (Krichhoff's law) at each node is satisfied for the given $\vect{x}(\mathcal{S})$.

    \item \emph{Step 3 (Check whether $(\mathcal{V}_D,\mathcal{S})\in\mathcal{F}$; computational complexity $O(m+n)$)}: Compute the undirected degree of each node in $(\mathcal{V}_D,\mathcal{S})$, i.e., sum of all the incident edges at $i$ regardless of direction. Insert all the nodes with undirected degree of one in queue $Q$, which was initialized empty. Repeat until $Q$ is empty: retrieve node $i$ from $Q$. Remove the sole incident edge $(i,j)$ of $i$ from $\mathcal{S}$ and  decrement the undirected degree of node $j$. If $j$'s degree becomes $1$ insert in $Q$. Confirm that $Q$ become empty after $|\mathcal{S}|$ element retrieval from $Q$.  
\end{itemize}
If all the aforementioned steps are checked out, the solution is confirmed and feasible. Since this verification steps complete in polynomial time with respect to the input size, the problem is in NP.

Next, we demonstrate weak NP-completeness by reduction from the Partition Problem, which is known to be weakly NP-complete \cite{CB-JC-BP:01}. The Partition Problem asks: given a multiset $\mathcal{A}=\{a_1 ,a_2, \cdots, a_n \}$ of positive integers, can $\mathcal{A}$ be divided into two subsets $\mathcal{A}_1$  and $\mathcal{A}_2$  such that the sum of elements in $\mathcal{A}_1$   equals the sum of elements in $\mathcal{A}_2$? (This implies that both sums must be equal to $\frac{1}{2}\sum_{i=1}^na_k$)~\cite{CB-JC-BP:01}. We proceed as follows:
\begin{enumerate}
    \item Create a graph $\mathcal{G}_D$ with $n + 3$ nodes: two source nodes $s_1, s_2$, and $n+1$ demand nodes $v_0,v_1, ..., v_n$.
    \item Add edges from both sources to all demand nodes, i.e., $(s_1, v_i)$ and $(s_2, v_i)$ for all $i\in\{0,1,\cdots, n\}$.
    \item Set the demand of each node $v_i$ to integer $d_i = a_i$, $i\in\{0,\cdots,n\}$.
    \item Set the edge capacity between node $s_i$, $i\in\{1,2\}$ and demand node $v_0$ to $a_0/2$. For the rest of demand nodes, set the edge capacities to $\bar{x}_{s_k,v_i} = a_i$, $k\in \{1,2\}$, $i\in\{1,\cdots,n\}$.
    \item Set the supply of each source node to $g_i=\frac{1}{2}\sum_{k=1}^na_k$, $i\in\{s_1,s_2\}$.
\end{enumerate}
Any feasible radial distribution network for this instance of problem~\eqref{eqn::problem1} must include a direct edge from each source $s_1$ and $s_2$ to node $v_0$ each supplying $a_0/2$ to meet the demand of this node. To guarantee a feasible radial distribution, the remaining sink nodes can only be exclusively connected to one of the sources (see Fig.~\ref{fig::np_complete}). Therefore, finding the feasible solution becomes an instance of Partition Problem. Since the Partition Problem is weakly NP-complete, problem~\eqref{eqn::problem1} is also weakly NP-hard. As it is also in NP, it is weakly NP-complete. 
\boxend\end{pf}

Creating a feasible solution for problem~\eqref{eqn::problem1} involves solving a combinatorial problem whose search space grows rapidly with network dimensionality. The next result provides the maximum number of possible combinations. Specifically, we demonstrate that each feasible distribution configuration, which is a polyforest, can be augmented by adding zero-flow edges to form an underlying undirected spanning tree of $\mathcal{G}_D$, and that the total number of such spanning trees (which includes those that do not yield a feasible flow) is given by the determinant of the cofactor of the Laplacian and defines the total number of possible feasible distribution~networks.

\begin{lem}\longthmtitle{Computational complexity of feasible solution of~\eqref{eqn::problem1}}\label{lem_size_solution}
    In a distribution network $\mathcal{G}_D$ with $n$ nodes, $m$ edges and $n_g$ source nodes the following assertions hold about the feasible solutions of problem~\eqref{eqn::problem1}: the number of feasible distribution configuration for problem~\eqref{eqn::problem1} is at most $\text{det}(\mathcal{L}^\prime(\mathcal{G}_D))$. Here, $\text{det}(\mathcal{L}^\prime(\mathcal{G}_D))$ is the determinant of the co-factor of the Laplacian of the graph $\mathcal{G}_D$.
\end{lem}
\begin{pf}
Let $\mathcal{G}_F = (\mathcal{V}_D, \mathcal{S})$ be a feasible distribution configuration, which is a solution to Problem~\eqref{eqn::problem1}. By Definition~\ref{def::radial}, $\mathcal{G}_F$ is a polyforest that includes all nodes $\mathcal{V}_D$ and has roots at all or some of $\mathcal{V}_g$\footnote{While source nodes in $\mathcal{V}_g$ are considered roots of the polyforest, it is possible for them to have incoming edges if their own supply $g_k$ is insufficient to satisfy all downstream demands they are responsible for, or if they act as relay nodes.}. The underlying undirected graph of $\mathcal{G}_F$, denoted $\mathcal{G}_{F,undir} = (\mathcal{V}_D, \mathcal{S}_{undir})$, is a spanning forest of $\mathcal{G}_D$ since $\mathcal{G}_F$ is acyclic and connects all nodes to the source set. If $\mathcal{G}_{F,undir}$ consists of $k > 1$ connected components, we can augment it by adding $k-1$ additional edges from the original network $\mathcal{E}_D \setminus \mathcal{S}_{undir}$ to connect these components. These additional edges must be chosen such that they do not form any cycles with the existing edges in $\mathcal{S}_{undir}$, thereby transforming $\mathcal{G}_{F,undir}$ into a single spanning tree of $\mathcal{G}_D$. For these newly added edges, we can assign a flow of zero ($x_{i,j}=0$). Therefore, any feasible distribution configuration can be represented by selecting edges from an underlying undirected spanning tree of $\mathcal{G}_D$.

The Matrix Tree Theorem (Kirchhoff's Theorem)~\cite{BB:98} states that for a connected undirected graph $\mathcal{G}_D$, the total number of its spanning trees is given precisely by the determinant of any cofactor of its Laplacian matrix, $\text{det}(\mathcal{L}^\prime(\mathcal{G}_D))$. Since every feasible distribution configuration can be augmented to form an underlying undirected spanning tree of $\mathcal{G}_D$, the number of feasible solutions is thus bounded and the number of feasible distribution configurations for problem~\eqref{eqn::problem1} is at most $\text{det}(\mathcal{L}^\prime(\mathcal{G}_D))$.\boxend
\end{pf}

A simple flow concept, starting from the pendant nodes and propagating demands/supplies inwards until no pendant edge remains, can be employed to find a feasible solution over a tree graph. This process can be completed in polynomial time, specifically $O(m+n)$, where $n$ is the number of nodes and $m$ is the number of edges. This approach is formalized in Algorithm~\ref{alg::tree_processor}, which guarantees either a feasible solution or a declaration of infeasibility. Furthermore, the subsequent result demonstrates that any existing solution on a tree graph is unique, and therefore, the solution generated by Algorithm~\ref{alg::tree_processor} is also unique.

\begin{algorithm}[t]
\small
\caption{\textsf{Tree-Processor}}
\begin{algorithmic}[1]
\Require Unidrected connected tree graph $\mathcal{G}(\mathcal{V},\mathcal{E})$, value vector $\vect{p}$
\State $\mathcal{S} \leftarrow \emptyset$
\While{exists pendant edges in $\mathcal{E}$}
\State Pick a pendant edge $(i,j)\in\mathcal{E}$, let $i$ denote the pendant node
\If{$\vect{p}_i>0$ and $|\vect{p}_i|\leq \bar{x}_{i,j}$}
\State $\mathcal{S}\leftarrow \mathcal{S}\cup\{(i\rightarrow j)\}$
\ElsIf{$\vect{p}_i=0$}
\State Do nothing \Comment{$(i,j)$ is not needed for flow transfer}
\ElsIf{$\vect{p}_i<0$ and $|\vect{p}_i|\leq \bar{x}_{i,j}$}
\State $\mathcal{S}\leftarrow\mathcal{S}\cup\{(j\rightarrow i)\}$
\Else
\State Terminate the entire function, and declare that no feasible solution exists
\EndIf
\State update $\vect{p}_j=\vect{p}_j+\vect{p}_i$
\State $\mathcal{V}\leftarrow \mathcal{V}\backslash\{i\}$ and $\mathcal{E}\leftarrow \mathcal{E}\setminus\{(i,j)\}$
\EndWhile
\State \Return $\mathcal{G}_F(\mathcal{V},\mathcal{S})$
\end{algorithmic}
\label{alg::tree_processor}
\vspace{0.1in}
\end{algorithm}

\begin{lem}\longthmtitle{Uniqueness of Flow and Radial Configuration in Trees}
\label{lem::uniqueness}
Given an undirected acyclic graph (a tree graph) where the nodes are either sink or source nodes with prespecified input and output values such that total input equals total output, the resulting flow distribution is unique. Consequently, the radial distribution configuration (polyforest) formed by these flows is also unique.
\end{lem}
\begin{pf}
For a connected acyclic graph with $N$ nodes, such as a tree graph, its incidence matrix $A$ has full column rank $(N-1)$~\cite{DH:22}. This implies that the null space of the system of flow conservation equations, $A \vect{x} = \vect{g} - \vect{d}$, is trivial. Given that the total input matches the total demand ($\sum_{i\in\mathcal{V}_g} g_i = \sum_{i\in\mathcal{V}_c} d_i$), the system is consistent. A consistent linear system with a trivial null space has a unique solution. Therefore, the flow vector $\vect{x}$ across all edges is uniquely determined.
Since the flow on each edge is uniquely determined, and a radial configuration (polyforest) is formed by the directed edges carrying non-zero flow, this flow distribution uniquely defines the set of active edges and their directions. Thus, the directed radial configuration is also unique.\boxend
\end{pf}

Building on Lemma~\ref{lem_size_solution} and Lemma~\ref{lem::uniqueness}, and utilizing Algorithm~\ref{alg::tree_processor}, which establishes that every feasible distribution configuration corresponds to an underlying undirected spanning tree of $\mathcal{G}_D$, a brute-force approach to Problem~\eqref{eqn::problem1} would involve enumerating all spanning trees of $\mathcal{G}_D$ and then checking each for feasibility. Due to the presence of capacity constraints, in the worst case, this exhaustive examination is theoretically necessary to guarantee finding a feasible solution. However, this is computationally intractable: the number of spanning trees can be exponential in $n$ (e.g., $n^{n-2}$ for a complete graph with $n$ nodes), leading to an exponential computational complexity for such a brute-force method and rendering it impractical for real-world distribution networks. To circumvent this challenge, we present an algorithm in the following section that efficiently identifies a feasible solution in polynomial time $\mathcal{O}(n^2 \log n)$.

\setlength{\textfloatsep}{3pt}
\begin{figure}[t]
  \centering
  \begin{subfigure}[t]{0.205\textwidth}
    \centering
    \raisebox{1.5em}{\includegraphics[width=\textwidth]{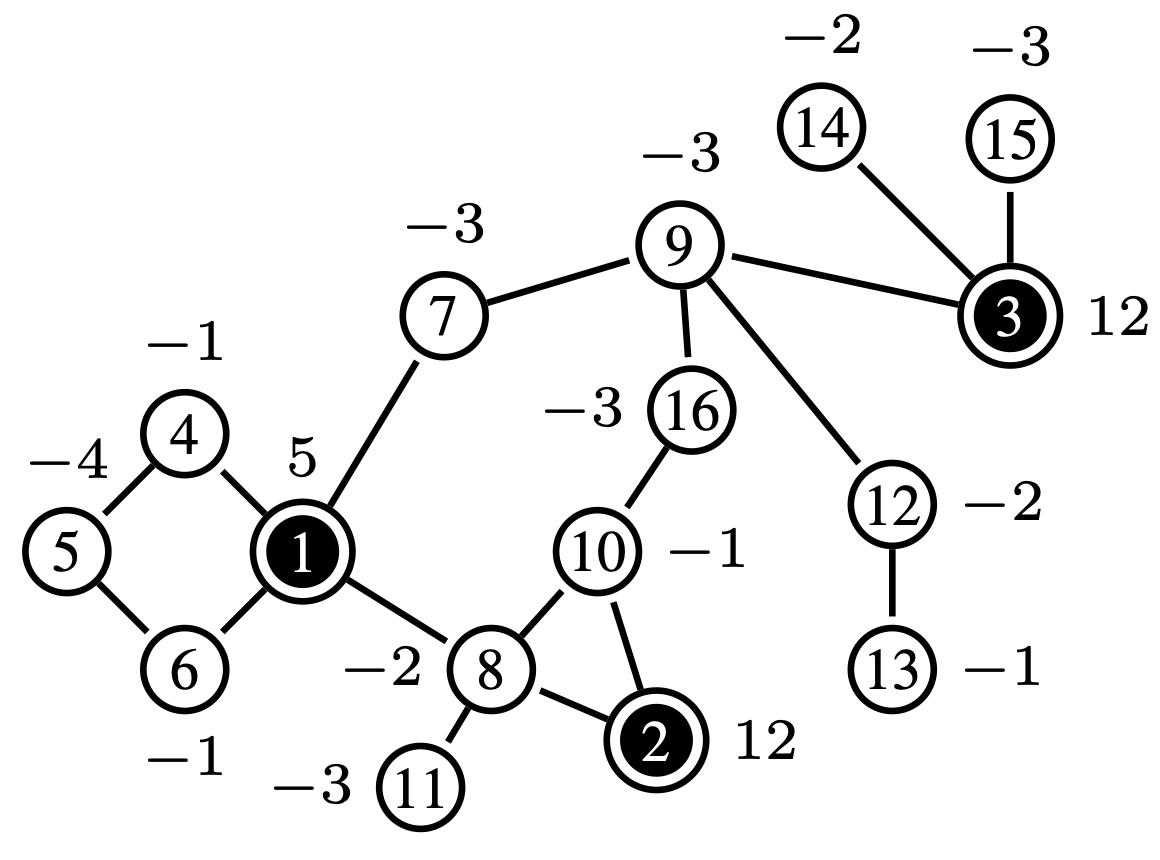}}
    \caption{{Distribution graph~$\mathcal{G}_D$}}
    \label{fig::original_graph}
  \end{subfigure}
  \begin{subfigure}[t]{0.26\textwidth}
    \centering
    \includegraphics[width=\textwidth]{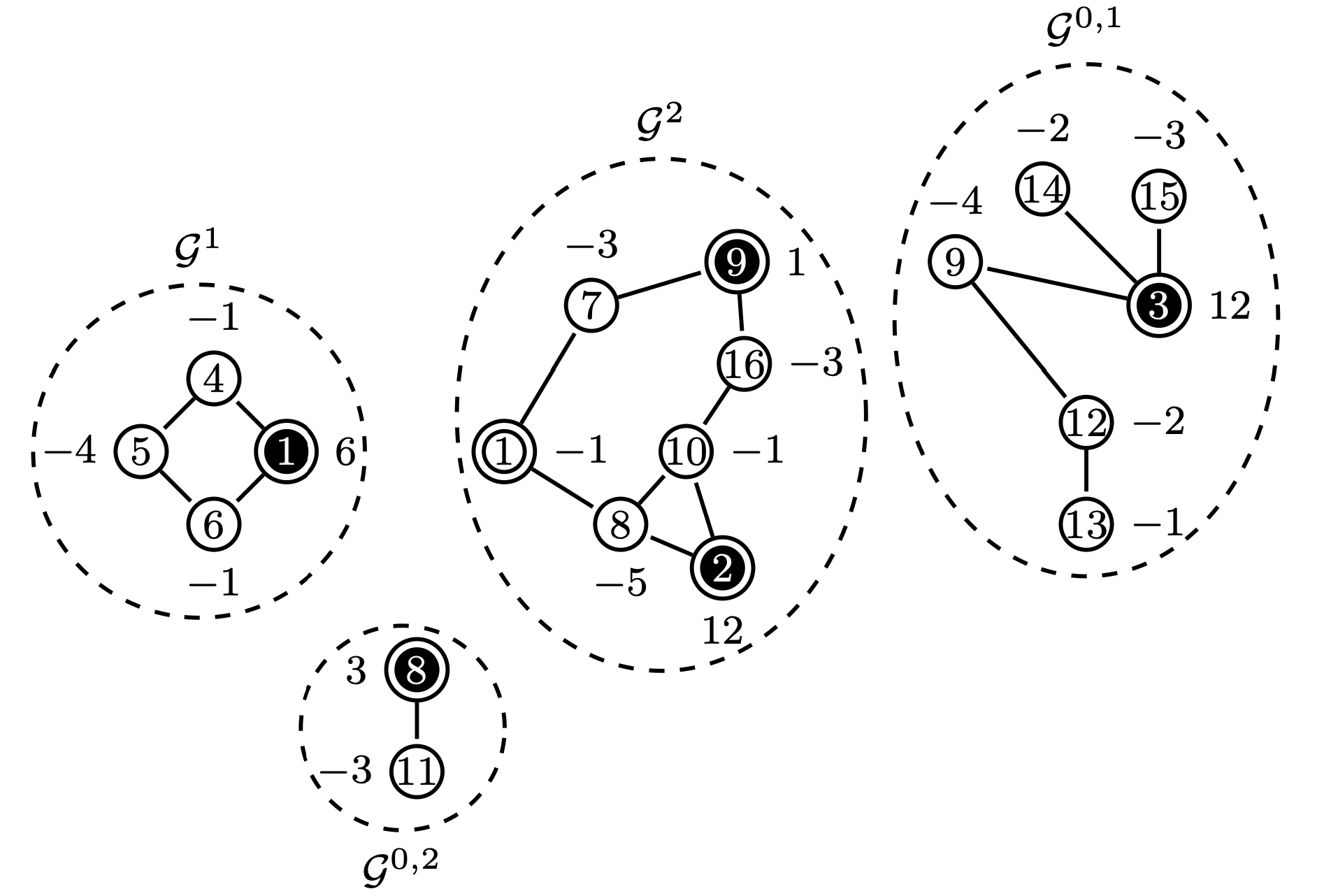}
    \caption{Partitioned subgraphs}
    \label{fig::extended_graph2}
  \end{subfigure}
  \caption{\footnotesize{A distribution network $\mathcal{G}_D$ (plot (a)) and partitioned subgraphs with adjusted nodal value at separation nodes (plot (b)).}}\vspace{0.1in}
  \label{fig::uniqueness}
\end{figure}
\section{Graph decomposition for parallel processing}
\label{sec::GraphDEcompose}
\vspace{-0.1in}
Distribution networks often contain nodes with low connectivity, creating what can be considered ``pendent" radial subgraphs that are intrinsic to the original network $\mathcal{G}_D$. These subgraphs must necessarily be part of any feasible solution, as they have no alternative paths for flow. Beyond these inherent structures, distribution networks frequently possess other articulation points (cut vertices) that allow us to further partition the graph into smaller connected components.

We propose a partitioning strategy that begins by first identifying the 2-core of the distribution graph $\mathcal{G}_D$, denoted by $\mathcal{G}_P$, and its leaf components. By partitioning the 2-core further at some of its articulation nodes (will be made more clear in proceeding section), we partition the entire network $\mathcal{G}_D$ into a set of disjoint connected components, $\mathcal{G}^{0,1},\cdots, \mathcal{G}^{0,K},\mathcal{G}^1,\ldots,\mathcal{G}^L$, where $\mathcal{G}^{0,i}$ are connected components of $\mathcal{G}^0$, the collection of all leaf components. The rest of $\mathcal{G}^j$s are derived from partitioning the 2-core subgraph $\mathcal{G}_P$ at some of its articulation points. An example of this decomposition is illustrated in Fig.~\ref{fig::uniqueness}. To ensure that a feasible flow can be found for each subproblem, the nodal values at the separation points must be adjusted such that each component becomes individually input/output balanced. This can be achieved through an iterative or recursive process, where the flow requirements of components with only a single separating node are first resolved and their net flow propagated to the adjacent component, similar to the \textsf{Tree-Processor} algorithm, until all subgraphs are balanced.

The next theorem establishes the theoretical foundation for this decomposition approach by proving a fundamental equivalence between the feasibility of the original problem and the feasibility of the balanced subproblems.

\begin{thm}\longthmtitle{Compositional Feasibility of Radial Distributions under Graph Partitioning}
\label{thm::partition}
Consider optimization problem~\eqref{eqn::problem1} under Assumption~\ref{assump::feas}. Let $\mathcal{G}_D$ be partitioned into disjoint components $\mathcal{G}^0, \mathcal{G}^1,\ldots,\mathcal{G}^L$ where $\mathcal{G}^0=\mathcal{G}_D\backslash\mathcal{G}_P$ and $\mathcal{G}^1,\ldots,\mathcal{G}^L$ are obtained by partitioning the 2-core $\mathcal{G}_P$ at $L\leq \bar{L}$ articulation nodes. After partitioning, nodal values are adjusted at separation nodes to ensure each component $\mathcal{G}^i$, $i\in\{0,\ldots,L\}$, is individually input/output balanced. Under this setting, a feasible radial distribution network exists for the original problem~\eqref{eqn::problem1} if and only if there exist feasible radial distributions for each balanced component $\mathcal{G}^\ell$, $\ell \in \{0, \ldots, L\}$.\boxend
\end{thm}
To conserve space, the proof is presented in Appendix~\ref{Appedix::Proofs}.

\noindent Beyond feasibility, the separable structure of the quadratic cost function across disjoint edge sets ensures that optimality is also preserved under this decomposition, as formalized in the following result.

\begin{cor}\longthmtitle{Optimality Preservation under Graph Partitioning}
\label{cor::optimality_partition}
Under the conditions of Theorem~\ref{thm::partition}, let $(\mathcal{S}^{\ell\star}, \vect{x}^{\ell\star})$ be an optimal solution to problem~\eqref{eqn::problem1} restricted to each balanced component $\mathcal{G}^\ell$, $\ell \in \{0, \ldots, L\}$. Then the combined solution $(\mathcal{S}^\star, \vect{x}^\star) = \left(\bigcup_{\ell=0}^L \mathcal{S}^{\ell\star}, \{\vect{x}^{\ell\star}\}_{\ell=0}^L\right)$ is an optimal solution to the original problem~\eqref{eqn::problem1} on $\mathcal{G}_D$.\boxend
\end{cor}
To conserve space, the proof is presented in Appendix~\ref{Appedix::Proofs}.

\vspace{0.5em}
\noindent Together, Theorem~\ref{thm::partition} and Corollary~\ref{cor::optimality_partition} provide the complete theoretical justification for decomposing large-scale radial distribution problems into smaller, independent subproblems that can be solved in parallel without loss of optimality.


\section{Algorithm Design for a Feasible Radial Distribution}
\label{sec::algorithm}
\vspace{-0.1in}
This section introduces our proposed \texttt{FORWARD} algorithm, designed to construct a radial configuration within a distribution network $\mathcal{G}_D$ that yields a feasible solution for the optimal distribution problem~\eqref{eqn::problem1}. \texttt{FORWARD} is built on the premise of  Theorem~\ref{thm::partition}, which establishes that a feasible solution for the entire network can be achieved by finding feasible solutions within its decomposed parts, provided these parts are made component-wise input-output balanced by adjusting values at their separation points. \texttt{FORWARD} comprises five main functions, which we will briefly introduce as we outline the algorithm's structure (Algorithm~\ref{alg::algorithm}), followed by a detailed explanation of each.

\begin{algorithm}[t]
\small
{\footnotesize
\caption{FORWARD}
\begin{algorithmic}[1]
\Require Bidirectional graph $\mathcal{G}_D$, value vector $\vect{p}$
\State $(\mathcal{S}^0,\mathcal{V}_g,\mathcal{G}_P,\vect{p}) \leftarrow$ \textsf{Pre-Processor}($\mathcal{G}_D,\vect{p}$)
\State $(\mathcal{G}^1,\mathcal{V}_g^1,\vect{p}^1),\cdots,(\mathcal{G}^L,\mathcal{V}_g^L,\vect{p}^L)\leftarrow$\textsf{Islander}($\mathcal{G}_P,\mathcal{V}_g,\vect{p}$)

\For{each partition $\ell\in\{1,\cdots,L\}$}
     \State $\mathbb{T}^\ell\gets \emptyset$ \Comment{Initializes the set of sampled polytrees}
    \For {$i\in\mathcal{V}^\ell_g$}\Comment{Initializes the sampled polytrees}
    \State $\mathcal{T}^\ell_i = \mathcal{G}(\{i\},\{\})$
    \State $\mathbb{T}^\ell\gets \mathbb{T}^\ell\cup \mathcal{T}^\ell_i$
    \EndFor
        \State $\mathcal{S}^\ell\leftarrow \emptyset$ \Comment{Initializes the sampled edges in $\mathcal{G}^\ell$}
        \State $e^\star\gets \{\}$
    \While {$|\mathcal{V}(\mathcal{S}^\ell)|\neq |\mathcal{V}^\ell|$} 
        \State ($\bar{\mathcal{G}}^\ell,\bar{\vect{p}}^\ell,\mathbb{T}^\ell)\leftarrow$\textsf{Net-Concad}($\mathcal{G}^\ell, \mathbb{T}^\ell,  \vect{p}^\ell,e^\star$)
        \State $e^\star\leftarrow$ \textsf{Sampler}($\bar{\mathcal{G}}^\ell,\vect{p}^\ell,\mathcal{S}^\ell$)
        \State $\mathcal{S}^\ell\leftarrow\mathcal{S}^\ell\cup\{e^\star\}$
    \EndWhile
    \If{at least one element of $\bar{\vect{p}}^\ell\neq 0$} 
       \State $\mathcal{S}^\ell \leftarrow$ \textsf{Rewire}($\bar{\mathcal{G}}^\ell,\vect{p}^\ell,\mathcal{S}^\ell$)
    \EndIf 
\EndFor
\State $\mathcal{S}\leftarrow \bigcup_{\ell=0}^L\mathcal{S}^\ell$ 
\State \Return Radial configuration $\mathcal{G}(\mathcal{V}(\mathcal{S}),\mathcal{S})$
\end{algorithmic}
\label{alg::algorithm}
}
\end{algorithm}

The \texttt{FORWARD} algorithm (Algorithm~\ref{alg::algorithm}) takes as input the distribution network graph $\mathcal{G}_D$, its initial nodal value vector $\vect{p}$, and its set of source nodes $\mathcal{V}_g$. The assumption is the problem~\eqref{eqn::problem1} is well-defined and Assumption~\ref{assump::feas} holds. The \texttt{FORWARD} algorithm first calls the \textsf{Pre-Processor} function (line 1 of Algorithm~\ref{alg::algorithm}).
\begin{itemize}
\item The \textsf{Pre-Processor} function simplifies the initial network by iteratively processing pendant nodes and redistributing their input/output to parent nodes. This reduces problem complexity by eliminating trivial components that must necessarily be included in any feasible solution.
\end{itemize}

Next, the \textsf{Islander} function is invoked (line 2 of Algorithm~\ref{alg::algorithm}).
\begin{itemize}
\item The \textsf{Islander} function partitions the 2-core subgraph $\mathcal{G}_P$ into disjoint subgraphs by splitting at articulation source nodes. This enables parallel processing and more efficient exploration of the solution space.
\end{itemize}

\texttt{FORWARD} then enters its parallel processing phase for each partitioned subgraph $\mathcal{G}^\ell$ to produce a feasible radial distribution configuration (lines 3-19). The algorithm initializes a set of polytrees $\mathbb{T}^\ell$ (lines 4-8), where each polytree $\mathcal{T}\in\mathbb{T}^\ell$ corresponds to a distinct source node from $\mathcal{V}_g^\ell$. These polytrees will incrementally `grow' through edge sampling by the \textsf{Sampler} and \textsf{Rewire} functions to eventually cover the entire subgraph. During this process, polytrees may merge when connecting edges are sampled. An empty sampled edge set $\mathcal{S}^\ell$ is also initialized for each subgraph $\mathcal{G}^\ell$ (line 9).

After the aforementioned steps, \texttt{FORWARD} calls the \textsf{Net-Concad} function.
\begin{itemize}
\item The \textsf{Net-Concad} function condenses each subgraph $\mathcal{G}^\ell$ by grouping polytrees into super source nodes and unsampled sink components into super sink nodes. This creates a quasi-bipartite structure $\bar{\mathcal{G}}^\ell$ that simplifies subsequent sampling and addresses the shortsightedness of the incremental greedy process.
\end{itemize}

After this initial call to \textsf{Net-Concad}, \texttt{FORWARD} enters an iterative loop (lines 11-15). Within each iteration, every dual graph $\bar{\mathcal{G}}^\ell$ must be updated due to polytree growth. Subsequently, the \textsf{Sampler} function is called.
\begin{itemize}
\item The \textsf{Sampler} selects an edge to add to the evolving radial configuration based on a set of complex prioritization criteria involving weights and network flow constraints to ensure a feasible and near-optimal process.
\end{itemize}

This iterative process, driven by the repeated application of \textsf{Net-Concad} and \textsf{Sampler}, dynamically modifies the super-node structure as polytrees grow or merge. These iterations continue until an attempt is made to cover all nodes, which requires at most $|\mathcal{V}^\ell|-1$ calls to the \textsf{Sampler} function.

When capacity constraints exist, the greedy selection process may prevent \textsf{Sampler} from achieving complete node coverage. In that case, the \textsf{Rewire} operator (lines 16-18) is activated.
\begin{itemize}
\item The \textsf{Rewire} function identifies edges that could provide missing flow to incompletely supplied nodes and detects which sampled edges are blocking their selection. Through swapping operations, it corrects previously sampled edges toward feasibility.
\end{itemize}

\setlength{\textfloatsep}{3pt}
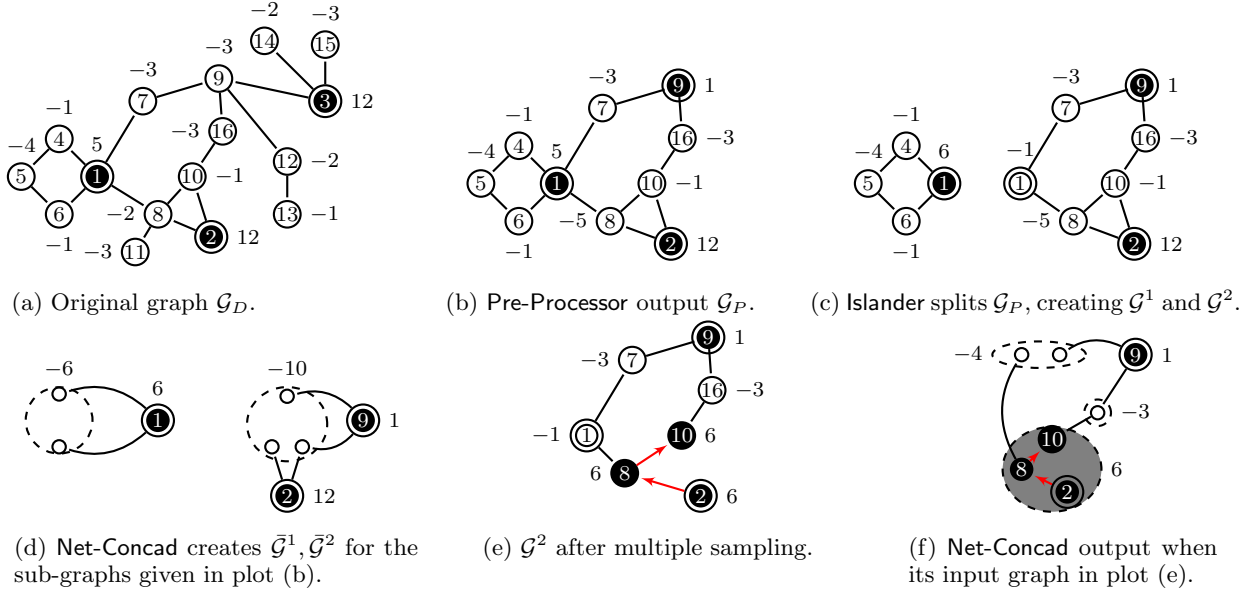
\begin{figure*}[t!]
    \centering  
    \scriptsize
    \begin{subfigure}[b]{0.2\textwidth}
        \centering
        \begin{tikzpicture}[>=stealth',shorten >=1pt,auto,node distance=1.5cm,
      thick,main node/.style={circle,fill=black,draw,minimum size=5pt,inner sep=0pt}]
          \node[main node, text=white, minimum size=8pt] (x1) at (0,0) {1};
          \node[main node, label=above:$5$, fill=none, minimum size=12pt] (x1111) at (0,0) {};
          \node[main node, text=white, minimum size=8pt] (x2) at (1.5,-0.8) {2};
          \node[main node, label=right:$12$, fill=none, minimum size=12pt] (x2222) at (1.5,-0.8) {};
          \node[main node, text=white, minimum size=8pt] (x3) at (3,1) {3};
          \node[main node, label=right:$12$, fill=none, minimum size=12pt, minimum size=12pt] (x3333) at (3,1) {};
          \node[main node, label=above:$-1$, fill=white, minimum size=10pt] (x4) at (-0.5,0.5) {4};
          \node[main node, label=above:$-4$, fill=white, minimum size=10pt] (x5) at (-1,0) {5};
          \node[main node, label=below:$-1$, fill=white, minimum size=10pt] (x6) at (-0.5,-0.5) {6};
          \node[main node, label=above:$-3$, fill=white, minimum size=10pt] (x7) at (0.6,1) {7};
          \node[main node, label=left:$-2$, fill=white, minimum size=10pt] (x8) at (0.8,-0.5) {8};
          \node[main node, label=right:$-1$, fill=white, minimum size=10pt] (x10) at (1.25,0) {10};
          \node[main node, label=left:$-3$, fill=white, minimum size=10pt] (x11) at (0.5,-1) {11};
          \node[main node, label=above:$-3$, fill=white, minimum size=10pt] (x9) at (1.6,1.3) {9};
          \node[main node, label=right:$-2$, fill=white, minimum size=10pt] (x12) at (2.5,0.2) {12};
          \node[main node, label=right:$-1$, fill=white, minimum size=10pt] (x13) at (2.5,-0.5) {13};
          \node[main node, label=above:$-3$, fill=white, minimum size=10pt] (x15) at (3,1.75) {15};
          \node[main node, label=above:$-2$, fill=white, minimum size=10pt] (x14) at (2.2,1.8) {14};
            \node[main node, label=left:$-3$, fill=white, minimum size=10pt] (x16) at (1.65,0.6) {16};
          \path[every node/.style={font=\sffamily\small}]
            (x1111) edge (x4)
            (x1111) edge (x6)
            (x1111) edge (x7)
            (x1111) edge (x8)
            (x4) edge (x5)
            (x5) edge (x6)
            (x7) edge (x9)
            (x8) edge (x11)
            (x8) edge (x10)
            (x9) edge (x16)
            (x16) edge (x10)
            (x2222) edge (x8)
            (x2222) edge (x10)
            (x3333) edge (x9)
            (x3333) edge (x14)
            (x3333) edge (x15)
            (x12) edge (x13)
            (x9) edge (x12);
        \end{tikzpicture}
         \caption{Original graph $\mathcal{G}_D$.}
        \label{fig:compact0}
    \end{subfigure}
    \hspace{15mm}
    \begin{subfigure}[b]{0.3\textwidth}
        \centering
        \begin{tikzpicture}[>=stealth',shorten >=1pt,auto,node distance=1.5cm,
      thick,main node/.style={circle,fill=black,draw,minimum size=5pt,inner sep=0pt}]
          \node[main node,text=white, minimum size=8pt] (x1) at (0,0) {1};
          \node[main node, label=above:$5$, fill=none, minimum size=12pt] (x1111) at (0,0) {};
          \node[main node, text=white, minimum size=8pt] (x2) at (1.5,-0.8) {2};
          \node[main node, fill=none, label=right:$12$, minimum size=12pt] (x2222) at (1.5,-0.8) {};
          \node[main node, label=above:$-1$, fill=white, minimum size=10pt] (x4) at (-0.5,0.5) {4};
          \node[main node, label=above:$-4$, fill=white, minimum size=10pt] (x5) at (-1,0) {5};
          \node[main node, label=below:$-1$, fill=white, minimum size=10pt] (x6) at (-0.5,-0.5) {6};
          \node[main node, label=above:$-3$, fill=white, minimum size=10pt] (x7) at (0.6,1) {7};
          \node[main node, label=left:$-5$, fill=white, minimum size=10pt] (x8) at (0.7,-0.5) {8};
          \node[main node, label=right:$-1$, fill=white, minimum size=10pt] (x10) at (1.25,0) {10};
       \node[main node,text=white, minimum size=8pt] (x9) at (1.6,1.3) {9};
          \node[main node, label=right:$1$, fill=none, minimum size=12pt] (x9999) at (1.6,1.3) {};
                \node[main node, label=right:$-3$, fill=white, minimum size=10pt] (x16) at (1.65,0.6) {16};
          \path[every node/.style={font=\sffamily\small}]
            (x1111) edge (x4)
            (x1111) edge (x6)
            (x1111) edge (x7)
            (x1111) edge (x8)
            (x4) edge (x5)
            (x5) edge (x6)
            (x7) edge (x9)
            (x8) edge (x10)
            (x9) edge (x16)
            (x16) edge (x10)
            (x2222) edge (x8)
            (x2222) edge (x10);
        \end{tikzpicture}
         \caption{\textsf{Pre-Processor} output $\mathcal{G}_P$.}
        \label{fig:compact1}
    \end{subfigure}
\begin{subfigure}[b]{0.32\textwidth}
        \centering
        \begin{tikzpicture}[>=stealth',shorten >=1pt,auto,node distance=1.5cm,
      thick,main node/.style={circle,fill=black,draw,minimum size=5pt,inner sep=0pt}]
          \node[main node,text=white, minimum size=8pt] (x1) at (0,0) {1};
          \node[main node, label=above:$6$, fill=none, minimum size=12pt] (x1111) at (0,0) {};
          \node[main node, text=black, fill=white, minimum size=8pt] (x1L) at (1,0) {1};
          \node[main node, label=above:$-1$, fill=none, minimum size=12pt] (x1111L) at (1,0) {};
          \node[main node, text=white, minimum size=8pt] (x2) at (2.5,-0.8) {2};
          \node[main node, fill=none, label=right:$12$, minimum size=12pt] (x2222) at (2.5,-0.8) {};
          \node[main node, label=above:$-1$, fill=white, minimum size=10pt] (x4) at (-0.5,0.5) {4};
          \node[main node, label=above:$-4$, fill=white, minimum size=10pt] (x5) at (-1,0) {5};
          \node[main node, label=below:$-1$, fill=white, minimum size=10pt] (x6) at (-0.5,-0.5) {6};
          \node[main node, label=above:$-3$, fill=white, minimum size=10pt] (x7) at (1.6,1) {7};
          \node[main node, label=left:$-5$, fill=white, minimum size=10pt] (x8) at (1.7,-0.5) {8};
          \node[main node, label=right:$-1$, fill=white, minimum size=10pt] (x10) at (2.25,0) {10};
       \node[main node,text=white, minimum size=8pt] (x9) at (2.6,1.3) {9};
          \node[main node, label=right:$1$, fill=none, minimum size=12pt] (x9999) at (2.6,1.3) {};
                \node[main node, label=right:$-3$, fill=white, minimum size=10pt] (x16) at (2.65,0.6) {16};
          \path[every node/.style={font=\sffamily\small}]
            (x1111) edge (x4)
            (x1111) edge (x6)
            (x1111L) edge (x7)
            (x1111L) edge (x8)
            (x4) edge (x5)
            (x5) edge (x6)
            (x7) edge (x9)
            (x8) edge (x10)
            (x9) edge (x16)
            (x16) edge (x10)
            (x2222) edge (x8)
            (x2222) edge (x10);
        \end{tikzpicture}
         \caption{\textsf{Islander} splits $\mathcal{G}_P$, creating $\mathcal{G}^1$ and~$\mathcal{G}^2$.}
        \label{fig:compact2}
    \end{subfigure}
    \begin{subfigure}[t]{0.3\textwidth}
        \centering
        \begin{tikzpicture}[>=stealth',shorten >=1pt,auto,node distance=1.5cm,
      thick,main node/.style={circle,fill=black,draw,minimum size=5pt,inner sep=0pt}]
          \node[main node, text=white, minimum size=8pt] (x11) at (-0.7,0) {1};
          \node[main node, label=above:$6$, fill=none, minimum size=12pt] (x111) at (-0.7,0) {};
          \node[main node, fill=white] (x2) at (-2,0.35) {};
          \node[main node, fill=white] (x4) at (-2,-0.35) {};
          \node[main node, fill=none, label=above:$-6$, dashed, minimum size=25pt] (x234) at (-2,0) {};
          \node[main node, fill=white] (x5) at (1,0.32) {};
          \node[main node, fill=white] (x6) at (0.8,-0.35) {};
          \node[main node, fill=white] (x8) at (1.2,-0.35) {};
          \node[main node, fill=none,label=above:$-10$, dashed, shape=ellipse, minimum width=30pt, minimum height=28pt] (x234) at (1,-0.05) {};
          \node[main node, minimum size=8pt, text=white] (x13) at (1,-1) {2};
          \node[main node, fill=none, minimum size=12pt, label=right:$12$] (x132) at (1,-1) {};
          \node[main node, minimum size=8pt, text=white] (x9) at (2,0) {9};
          \node[main node, fill=none, minimum size=12pt, label=right:$1$] (x92) at (2,0) {};
          \path[every node/.style={font=\sffamily\small}]
            (x111) edge[bend right] (x2)
            (x111) edge[bend left] (x4)
            (x132) edge (x6)
            (x132) edge (x8)
            (x92) edge[bend right] (x5)
            (x92) edge[bend left] (x8);
        \end{tikzpicture}
        \caption{{\textsf{Net-Concad} creates $\bar{\mathcal{G}}^1,\bar{\mathcal{G}}^2$ for the sub-graphs given in plot (b).}}
        \label{fig:compact3}
    \end{subfigure}
        \hspace{6mm}
      \begin{subfigure}[t]{0.25\textwidth}
        \centering
        \begin{tikzpicture}[>=stealth',shorten >=1pt,auto,node distance=1.5cm,
      thick,main node/.style={circle,fill=black,draw,minimum size=5pt,inner sep=0pt}]
          \node[main node,fill=none, minimum size=8pt] (x1) at (0,0) {1};
          \node[main node, label=left:$-1$, fill=none, minimum size=12pt] (x1111) at (0,0) {};
          \node[main node, text=white, minimum size=8pt] (x2) at (1.5,-0.8) {2};
          \node[main node, fill=none, label=right:$6$, minimum size=12pt] (x2222) at (1.5,-0.8) {};
          \node[main node, label=left:$-3$, fill=white, minimum size=10pt] (x7) at (0.6,1) {7};
          \node[main node, text=white, minimum size=10pt, label=left:$6$] (x8) at (0.5,-0.5) {8};
          \node[main node, text=white, minimum size=10pt, label=right:$6$] (x10) at (1.25,0) {10};
       \node[ main node,  text=white, minimum size=8pt] (x9) at (1.6,1.3) {9};
          \node[main node, label=right:$1$, fill=none, minimum size=12pt] (x9999) at (1.6,1.3) {};
                \node[main node, label=right:$-3$, fill=white, minimum size=10pt] (x16) at (1.65,0.6) {16};
          \path[every node/.style={font=\sffamily\small}]
            (x1111) edge (x7)
            (x1111) edge (x8)
            (x7) edge (x9)
            (x9) edge (x16)
            (x16) edge (x10);
            \draw[arrow_e,red] (x8) to (x10);
            \draw[arrow_e,red] (x2222) to (x8);
        \end{tikzpicture}
         \caption{{$\mathcal{G}^2$ after multiple sampling.}}
       \label{fig:compact4}
    \end{subfigure}
        \hspace{10mm}
    \begin{subfigure}[t]{0.23\textwidth}
        \centering
        \begin{tikzpicture}[>=stealth',shorten >=1pt,auto,node distance=1.5cm,
      thick,main node/.style={circle,fill=black,draw,minimum size=5pt,inner sep=0pt}]
          \node[main node, fill=white] (x5) at (1,0) {};
          \node[main node, fill=white] (x6) at (0.5,0) {};
          \node[main node, fill=none,label=left:$-4$, dashed, shape=ellipse, minimum width=35pt, minimum height=10pt] (x234) at (0.73,0) {};
          \node[main node, fill=none,label=right:$6$, dashed, shape=ellipse, minimum width=37pt, minimum height=32pt, fill=gray] (x882) at (0.9,-1.52) {};
          \node[main node, fill=white] (x8) at (1.5,-0.77) {};
          \node[main node, fill=none,label=right:$-3$, dashed, shape=ellipse, minimum width=10pt, minimum height=10pt] (x888) at (1.5,-0.77) {};
          \node[main node, minimum size=8pt, text=white] (x13) at (1.1,-1.82) {2};
          \node[main node, fill=none, minimum size=12pt] (x132) at (1.1,-1.82) {};
          \node[main node, minimum size=8pt, text=white] (x88) at (0.5,-1.52) {8}; 
          \node[main node, minimum size=8pt, text=white] (x10) at (0.9,-1.12) {10}; 
          \node[main node, minimum size=8pt, text=white] (x9) at (2,0) {9};
          \node[main node, fill=none, minimum size=12pt, label=right:$1$] (x92) at (2,0) {};
          \path[every node/.style={font=\sffamily\small}]
            (x88) edge[arrow_e,red] (x10)
            (x8) edge (x10)
            (x88) edge[bend left] (x6)
            (x92) edge[bend right] (x5)
            (x92) edge (x8);
             \tikzset{arrow_e/.style = {->,> = latex'}}
            \draw[arrow_e,red] (x132) to (x88);
        \end{tikzpicture}
    \caption{{\textsf{Net-Concad} output when its input graph in plot (e).}}
       \label{fig:compact5}
    \end{subfigure}
    \hfill
    \caption{{\footnotesize Demonstration of how \textsf{Pre-Processor},\textsf{Islander} and \textsf{Net-Concad} function. Filled solid nodes represent sources while others represent the sinks. The number next to the nodes show corresponding $p_i$. In plot (b) $p_9$ and $p_8$ are adjusted to reflect, respectfully, excess input to deliver through node 9 to $\mathcal{G}_p$ and extra demand at node $8$ to supply to the removed pendent node $11$, being {\footnotesize $\mathcal{S}=\{(9\rightarrow 12),(12\rightarrow 13),(8\rightarrow 11),(3\rightarrow 14),(3\rightarrow 15)\}$}. In plot (c), after partitioning the graph, node $1$ assumes different roles in each sub-graph. In plot (d) and (f) dashed circles represent super nodes. In plot (e), $\mathcal{T}_1=\mathcal{G}(\{9\},\{\})$ and $\mathcal{T}_2=\mathcal{G}(\{2,8,10\},(2\to8),(8\to10))$; note that here, \textsf{Sampler} has removed edge $(2,10)$ to avoid cycle. Being {\footnotesize $\mathcal{S}=\{(9\rightarrow 12),(12\rightarrow 13),(8\rightarrow 11),(3\rightarrow 14),(3\rightarrow 15),(2\rightarrow 8),(8\rightarrow 10)\}$. 
    }}}    
    \label{fig:connected_comps}
\end{figure*}

\begin{algorithm}[t]
\small
{\footnotesize
\caption{\textsf{Pre-Processor}}
\begin{algorithmic}[1]
\Require 
Unidrected graph $\mathcal{G}_D$, value vector $\vect{p}$
\State $\mathcal{S} \leftarrow \emptyset$
\State $\mathcal{E}_P \leftarrow \mathcal{E}_D$
\While{exists pendant nodes}
\State Pick a pendent edge $(i,j)\in\mathcal{E}_P$, let $i$ denote the pendent node
\If{$\vect{p}_i>0$}
\State $\mathcal{S}\leftarrow \mathcal{S}\cup\{(i\rightarrow j)\}$ 
\ElsIf{$\vect{p}_i=0$}
\State Do nothing \Comment{$(i,j)$ is not needed for flow transfer}
\Else
\State $\mathcal{S}\leftarrow\mathcal{S}\cup\{(j\rightarrow i)\}$ 
\EndIf
\State update $\vect{p}_j=\vect{p}_j+\vect{p}_i$
\State $\mathcal{E}_P\leftarrow \mathcal{E}_P\setminus\{(i,j)\}$
\EndWhile
\State $\mathcal{G}_P\leftarrow \mathcal{G}(\mathcal{V}(\mathcal{E}_P),\mathcal{E}_P)$
\State Reshape $\vect{p}$ to include only the elements corresponding to nodes of $\mathcal{G}_P$
\State 
$\mathcal{V}_g=\{i\in\mathcal{V}(\mathcal{E}_p)|\vect{p}_i>0\}$ \Comment{Reshape $\mathcal{V}_g$ according to $\mathcal{G}_P$}
\State \Return $\mathcal{S}$, $\mathcal{G}_P$, $\vect{p}$, $\mathcal{V}_g$
\end{algorithmic}
\label{alg::preprocessor}
}
\end{algorithm}

\subsection{Algorithm's components}
\vspace{-0.1in}
\subsubsection{Pre-Processor}
\vspace{-0.1in}
The \textsf{Pre-Processor} function (Algorithm~\ref{alg::preprocessor}) takes the original distribution network graph $\mathcal{G}_D$ and the associated initial nodal value vector $\vect{p}$ as inputs. The function identifies existing radial subgraphs in $\mathcal{G}_D$ and constructs a feasible radial distribution over these subgraphs by systematically processing trivially included components.

Specifically, the \textsf{Pre-Processor} samples every edge connected to a pendant node (a node with degree one), such as edges $(8,11)$, $(9,12)$, $(12,13)$, $(3,14)$, $(3,15)$, and $(9,3)$ in Fig.~\ref{fig:compact0}. For each pendant node $i$ connected to node $j$, the algorithm adds the appropriately directed edge to the solution set $\mathcal{S}$ based on the sign of $\vect{p}_i$ (lines 5-11); Notice that if $\vect{p}_i = 0$, no edge is needed for flow transfer. The pendant node is then removed from the graph, and its input/output value is redistributed to its parent node (line 12). This process iterates until no pendant nodes remain in the updated graph. The result of applying the \textsf{Pre-Processor} function to the network in Fig.~\ref{fig:compact0} is depicted in Fig.~\ref{fig:compact1}. The \textsf{Pre-Processor} function  is similar to Algorithm~\ref{alg::tree_processor}, but due to Assumption~\ref{assump::feas} that guarantees the existence of a feasible solution and Theorem~\ref{thm::hardness} that establishes feasibility within radial components of 
$\mathcal{G}_D$, there is no need to check whether capacity constraints are satisfied during this process. 

After \textsf{Pre-Processor} simplification, all nodes in the resulting subgraph $\mathcal{G}_P$ have a degree of at least two, making $\mathcal{G}_P$ the 2-core subgraph of $\mathcal{G}_D$. The function returns the sampled edge set $\mathcal{S}$, the 2-core subgraph $\mathcal{G}_P$, the updated nodal value vector $\vect{p}$, and the updated source node set $\mathcal{V}_g$ (lines 15-18).

\subsubsection{Islander}
\label{sec:islander}
\vspace{-0.1in}
The \textsf{Islander} function (Algorithm~\ref{alg::islander}) takes as input the 2-core subgraph $\mathcal{G}_P$ generated by \textsf{Pre-Processor}, along with the corresponding source nodes $\mathcal{V}_g$ and nodal value vector $\vect{p}$. Note that $\mathcal{V}_g$ may differ from the original source set of $\mathcal{G}_D$ due to the \textsf{Pre-Processor}'s redistribution of nodal values (line 17 of Algorithm~\ref{alg::preprocessor}).

\textsf{Islander} partitions the the 2-core subgraph $\mathcal{G}_P$ into disjoint subgraphs $\mathcal{G}^\ell$, $\ell\in\{1,\ldots,L\}$, by splitting at articulation source nodes ($L\leq |\mathcal{V}_g|$, since not all source nodes are articulation nodes). The choice to partition only at source articulation nodes, rather than sink articulation nodes, reduces computational complexity while maintaining effectiveness, as sink articulations are naturally handled by the subsequent \textsf{Net-Concad} operation.

After partitioning, the value vector $\vect{p}$ of each articulation node is adjusted to balance the overall input-output in each subgraph through a simple LP problem or a process explained in Section~\ref{sec::GraphDEcompose}. This process can change the role of articulation nodes from source to sink or vice versa within different subgraphs. For example, applying \textsf{Islander} to the graph in Fig.~\ref{fig:compact1} (with source nodes $1$, $2$, and $9$) creates two subgraphs as shown in Fig.~\ref{fig:compact2}, where node $1$ has different roles in each subgraph.

{\small
\begin{algorithm}[t]
{\footnotesize
\caption{\textsf{Islander}}
\begin{algorithmic}[1]
\Require The 2-core subgraph $\mathcal{G}_P$, source nodes $\mathcal{V}_g$, value vector $\vect{p}$
\State Find nodes in $\mathcal{V}_g$ that are articulation nodes of $\mathcal{G}_P$
\State Partition $\mathcal{G}_P$ at articulation nodes creating $\mathcal{G}^\ell=\mathcal{G}(\mathcal{V}^\ell,\mathcal{E}^\ell)$, $\ell\in\{1,\cdots,L\}$
\State Define an value vector $\vect{p}^\ell$, which contains the elements of $\vect{p}$ that corresponds to nodes in $\mathcal{G}^\ell$, $\ell\in\{1,\cdots,L\}$
\State Compute and update $\vect{p}_i^\ell$ for each articulation point $i$ used in partitioning $\mathcal{G}_P$ in each partition $\ell$
\State $\mathcal{V}_g^\ell=\{i\in\mathcal{V}^\ell|\vect{p}_i^\ell>0\}$, ~~$\ell\in\{1,\cdots,L\}$
\State \Return $(\mathcal{G}^1,\mathcal{V}_g^1,\vect{p}^1),\cdots,(\mathcal{G}^L,\mathcal{V}_g^L,\vect{p}^L)$ 
\end{algorithmic}
\label{alg::islander}
}
\end{algorithm}
}

\subsubsection{Net-Concad}
\label{sec:net-concad}
\vspace{-0.1in}
The \textsf{Net-Concad} function (Algorithm~\ref{alg::netconcad}) takes as input a graph $\mathcal{G}^\ell$ and the set of existing sampled radial polytrees $\mathbb{T}^\ell$ within $\mathcal{G}^\ell$ and a newly sampled edge $i\to j$ by the \textsf{Sampler} function, which will be explained next. 
Initially, each $\mathcal{T}^\ell_i\in\mathbb{T}^\ell$ is a distinct source node from $\mathcal{V}_g^\ell$. As new edges are sampled by the \textsf{Sampler}, the polytrees incrementally expand to cover the entire $\mathcal{G}^\ell$. In the process, some or all of the polytrees can merge and form larger polytrees. 

 \textsf{Net-Concad} carries two functions: one is to update the exiting polytree set after an edge is sampled by the \textsf{Sampler}, which is done by expanding the correspondent polytree in $\mathbb{T}^\ell$ or merging them if the edge is between to existing sampled polytrees. Note that in the absence of the capacity bound, every node within a polytree inherently can have access to the remaining positive input in the source nodes of the polytree to supply to its neighboring sink nodes (alternatively becoming into a ``super sink" if a new node is added that demands more output than the polytree's current total nodal value).

 The second mission of \textsf{Net-Concad} function is to condense each subgraph $\mathcal{G}^\ell$ into a dual graph $\bar{\mathcal{G}}^\ell$ featuring cluster nodes. This condensation is achieved by grouping nodes within each formed polytree into a single super `sampled' node, and by grouping each connected sink components formed after removing the super `sampled' nodes into super un-sampled nodes, using \textsf{Conncomp} method described in~\cite{conncomp} (as depicted in Figures~\ref{fig:compact3} and~\ref{fig:compact5}). This new dual graph $\bar{\mathcal{G}}^\ell$ is quasi-bipartite, with super `sampled' nodes forming one partition and super `un-sampled' nodes forming the other. Formally, we define the dual concatenated  graph as $\bar{\mathcal{G}}^\ell(\bar{\mathcal{V}}^\ell,\bar{\mathcal{E}}^\ell)$, where $\bar{\mathcal{V}}^\ell$ is the node set consisted of sample and un-sampled super nodes and  $\bar{\mathcal{E}}^\ell$ is the edge set where every edge $e\in \bar{\mathcal{E}}^\ell$ is a tuple $(u', v', s, t)$, where $u'$ and $v'$ are supernodes in $\bar{\mathcal{V}}$, $s$ is a subnode in the supernode represented by $u'$ and $t$ is a subnode in the supernode $v'$. $\bar{\mathcal{G}}^\ell$ is an undirected graph; when an edge $(u', v', s, t)$ is between a super sampled node and a super un-sampled node, we let $u'$ be the super sampled node and $v'$ be the super un-sampled node. We define the \emph{collective nodal value} of any \emph{un-sampled} super node $v$ as $\bar{\vect{p}}^\ell_v=\sum_{i\in \mathcal{W}^\ell_v}\vect{p}_i^\ell$, where $\mathcal{W}^\ell_v$ is the set of nodes in super node $v$. We set $\bar{\vect{p}}^\ell_u=0$ for any any \emph{sampled} super node $u$. 

{\small
\begin{algorithm}[t]
{\footnotesize
\caption{\textsf{Net-Concad}}
\begin{algorithmic}[1]
\Require sub-graph $\mathcal{G}^\ell$, set of existing polytrees $\mathbb{T}^\ell$, nodal value vector $\vect{p}^\ell$ and newly sampled edge $i\to j$.
\State $\mathbb{T}^\ell \gets\textsf{Tree-Update}(\mathbb{T}^\ell,i\to j)$ \Comment{See Appendix~\ref{Appedix::subroutines}}
\State $\vect{p}^\ell\gets$ \textsf{Update}($\vect{p}^\ell,i\rightarrow j$) 
\State Define $\mathcal{V}^\ell_s\leftarrow\{i: \forall i\in\mathcal{V}(\mathbb{T}^\ell) \}$ and $\mathcal{V}^\ell_u\leftarrow\{i: \forall i\in\mathcal{V}^\ell\setminus\mathcal{V}(\mathbb{T}^\ell)\}$
\State Define $\mathcal{G}^\ell_s\leftarrow\mathcal{G}(\mathcal{V}^\ell_s,\mathcal{E}(\mathcal{V}^\ell_s))$ and $\mathcal{G}^\ell_u\leftarrow\mathcal{G}(\mathcal{V}^\ell_u,\mathcal{E}(\mathcal{V}_u))$ 
\State $\bar{\mathcal{G}}^\ell_s,\bar{\mathcal{G}}^\ell_u \leftarrow$ Run \textsf{Conncomp} over $\mathcal{G}^\ell_s$ and $\mathcal{G}^\ell_u$
\State  $\bar{\mathcal{G}}^\ell\leftarrow\mathcal{G}(\bar{\mathcal{V}}^\ell_s\cup\bar{\mathcal{V}}^\ell_u, \mathcal{E}^\ell_P\setminus\{\mathcal{E}(\mathcal{V}^\ell_s\cup\mathcal{V}^\ell_u)\})$\Comment{Re-connect sub-graphs}
\State \Return $\bar{\mathcal{G}}^\ell, \vect{p}^\ell$
\end{algorithmic}
\label{alg::netconcad}
}
\end{algorithm}
}


\subsubsection{Sampler}
\vspace{-0.1in}
The \textsf{Sampler} function (Algorithm~\ref{alg::sampler}) takes as input the graph $\bar{\mathcal{G}}^\ell$, the current nodal value vector $\vect{p}^\ell$ and the grown polytrees, 
and returns a newly sampled edge to be added to one of the polytrees.  \textsf{Sampler} employs $\textsf{Weight}(\bar{\mathcal{G}}^\ell,\bar{\mathcal{G}}^\ell,\vect{p}^\ell, \bar{\vect{p}}^\ell)$ function to determine the ``most probable" edge (edge with the highest weight) to include from a set of candidate edges. 

The \textsf{Weight} function assigns a weight $w_{i,j}>0$ to every edge $(u,v,i,j)\in\bar{\mathcal{E}}_\ell$ if only $\bar{\vect{p}}_u^\ell>0$ or only $\bar{\vect{p}}_v^\ell>0$, i.e., when $\textsf{XOR}(\bar{\vect{p}}_u^\ell>0, \bar{\vect{p}}_v^\ell>0)=1$; otherwise $w_{i,j}=0$. The weighting process assigns a positive weight to an edge if and only if it connects a super-sampled node with a positive collective nodal value to another super-node (sampled or un-sampled) with a non-positive collective nodal value. The \textsf{Weight} function returns all the tuples $((u,v,i,j), w_{i,j})$ with positive weight. This ensures the \textsf{Sampler} only considers edges where one end node consistently has a positive input to supply, thereby pushing flow from sources to sinks. These weights dynamically change after each sampling event, being defined by the specific problem context and the physical characteristics of the flow. Section~\ref{sec::numerical} provides a concrete example of the electric power distribution \textsf{Weight} function.

Finally, when \textsf{Sampler} is called in \texttt{FORWARD} it will sample an edge by strategically using a priority queue. It first populates this queue with tuples $((u,v,i,j), w_{i,j})$ outputted from $\textsf{Weight}$. This queue is primarily sorted by descending weight, but further prioritized: highest for $v$ being a pendent super un-sampled node or $(u,v)$ that possesses sufficient capacity to carry the required flow. Meanwhile \textsf{Weight} function definition advocates for the optimality on the sampled solution, the priority queue is crucial for yielding into a feasible configuration. 

\begin{figure}[t]
    \centering
        \begin{subfigure}{0.1\textwidth}
           \centering
            \begin{tikzpicture}[>=stealth',shorten      >=1pt,auto,node distance=1.5cm,
            thick,main node/.style={circle,fill=black,draw,minimum size=5pt,inner sep=0pt}]
              \node[main node, fill=black, text=white, minimum size=6pt] (x1) at (0,0) {\scriptsize 1};
              \node[main node, fill=none, label=above:{\tiny$3$}, minimum size=9pt] (x1111) at (0,0) {};
              \node[main node, fill=none, dashed, shape=rectangle, minimum width=15pt, minimum height=12pt] (x11) at (0,0) {};
              \node[main node, fill=white, label=above:{\tiny$-15$}, minimum size=8pt] (x2) at (0.7,0) {\scriptsize 2};
              \node[main node, fill=black, text=white, minimum size=6pt] (x3) at (1.4,0) {\scriptsize 3};
              \node[main node, fill=none, label=above:{\tiny$10$}, minimum size=9pt] (x3333) at (1.4,0) {};
              \node[main node, fill=none, dashed, shape=rectangle, minimum width=15pt, minimum height=12pt] (x33) at (1.4,0) {};
              \node[main node, fill=white, label=below:{\tiny$-11$}, minimum size=8pt] (x4) at (0,-1) {\scriptsize 4};
              \node[main node, fill=black, text=white, minimum size=6pt] (x5) at (0.7,-1) {\scriptsize 5};
              \node[main node, fill=none, label=below:{\tiny$20$}, minimum size=9pt] (x5555) at (0.7,-1) {};
              \node[main node, fill=none, dashed, shape=rectangle, minimum width=15pt, minimum height=12pt] (x55) at (0.7,-1) {};
              \node[main node, fill=white, label=below:{\tiny$-7$}, minimum size=8pt] (x6) at (1.4,-1) {\scriptsize 6};
                \draw (x1111) -- (x2);
                \draw (x2) -- (x3333);
                \draw (x1111) -- (x4);
                \draw (x4) -- (x5555);
                \draw (x5555) -- (x6);
                \draw (x6) -- (x3333);
            \end{tikzpicture}
            \caption{{ ${\mathcal{G}}^\ell$}}
        \end{subfigure}
    \hspace{5mm}
        \begin{subfigure}{0.1\textwidth}
        \centering
            \begin{tikzpicture}[>=stealth',shorten >=1pt,auto,node distance=1.5cm,
            thick,main node/.style={circle,fill=black,draw,minimum size=5pt,inner sep=0pt}]
              \node[main node, fill=black, text=white, minimum size=6pt] (x1) at (0,0) {\scriptsize 1};
              \node[main node, fill=none, label=above:{\tiny$3$}, minimum size=9pt] (x1111) at (0,0) {};
              \node[main node, fill=none, dashed, shape=rectangle, minimum width=15pt, minimum height=12pt] (x11) at (0,0) {};
              \node[main node, fill=white, label=above:{\tiny$-15$}, minimum size=8pt] (x2) at (0.7,0) {\scriptsize 2};
              \node[main node, fill=black, text=white, minimum size=6pt] (x3) at (1.4,0) {\scriptsize 3};
              \node[main node, fill=none, label=above:{\tiny$10$}, minimum size=9pt] (x3333) at (1.4,0) {};
              \node[main node, fill=none, dashed, shape=rectangle, minimum width=15pt, minimum height=12pt] (x33) at (1.4,0) {};
              \node[main node, fill=black, text=white, label=below:{\tiny$9$}, minimum size=8pt] (x4) at (0,-1) {\scriptsize 4};
              \node[main node, fill=black, text=white, minimum size=6pt] (x5) at (0.7,-1) {\scriptsize 5};
              \node[main node, fill=none, dashed, shape=rectangle, minimum width=32pt, minimum height=12pt] (x45) at (0.35,-1) {};
              \node[main node, fill=none, label=below:{\tiny$9$}, minimum size=9pt] (x5555) at (0.7,-1) {};
              \node[main node, fill=white, label=below:{\tiny$-7$}, minimum size=8pt] (x6) at (1.4,-1) {\scriptsize 6};
                \draw (x1111) -- (x2);
                \draw (x2) -- (x3333);
                \draw (x1111) -- (x4);
                \draw[->,red] (x5555) -- (x4);
                \draw (x5555) -- (x6);
                \draw (x6) -- (x3333);
            \end{tikzpicture}
            \caption{{ $1^{\textup{st}}$ step}}
        \end{subfigure}
        \hspace{5mm}
        \begin{subfigure}{0.1\textwidth}
        \centering
            \begin{tikzpicture}[>=stealth',shorten >=1pt,auto,node distance=1.5cm,
            thick,main node/.style={circle,fill=black,draw,minimum size=5pt,inner sep=0pt}]
              \node[main node, fill=black, text=white, minimum size=6pt] (x1) at (0,0) {\scriptsize 1};
              \node[main node, fill=none, label=above:{\tiny$3$}, minimum size=9pt] (x1111) at (0,0) {};
              \node[main node, fill=none, dashed, shape=rectangle, minimum width=15pt, minimum height=12pt] (x11) at (0,0) {};
              \node[main node, fill=white, label=above:
              {\tiny$-15$}, minimum size=8pt] (x2) at (0.7,0) {\scriptsize 2};
              \node[main node, fill=black, text=white, minimum size=6pt] (x3) at (1.4,0) {\scriptsize 3};
              \node[main node, fill=none, label=above:{\tiny$10$}, minimum size=9pt] (x3333) at (1.4,0) {};
              \node[main node, fill=none, dashed, shape=rectangle, minimum width=15pt, minimum height=12pt] (x33) at (1.4,0) {};
              \node[main node, fill=black, text=white, label=below:{\tiny$2$}, minimum size=8pt] (x4) at (0,-1) {\scriptsize 4};
              \node[main node, fill=black, text=white, minimum size=6pt] (x5) at (0.7,-1) {\scriptsize 5};
              \node[main node, fill=none, dashed, shape=rectangle, minimum width=52pt, minimum height=12pt] (x45) at (0.7,-1) {};
              \node[main node, fill=none, label=below:{\tiny$2$}, minimum size=9pt] (x5555) at (0.7,-1) {};
              \node[main node, fill=black, text=white, label=below:{\tiny$2$}, minimum size=8pt] (x6) at (1.4,-1) {\scriptsize 6};
                \draw (x1111) -- (x2);
                \draw (x2) -- (x3333);
                \draw (x1111) -- (x4);
                \draw[->,red] (x5555) -- (x4);
                \draw[->,red] (x5555) -- (x6);
                \draw (x6) -- (x3333);
            \end{tikzpicture}
            \caption{{$2^{\textup{nd}}$ step}}
        \end{subfigure}
        \begin{subfigure}{0.1\textwidth}
        \centering
            \begin{tikzpicture}[>=stealth',shorten >=1pt,auto,node distance=1.5cm,
            thick,main node/.style={circle,fill=black,draw,minimum size=5pt,inner sep=0pt}]
              \node[main node, fill=white, text=black, minimum size=6pt] (x1) at (0,0) {\scriptsize 1};
              \node[main node, fill=none, label=above:{\tiny$-12$}, minimum size=9pt] (x1111) at (0,0) {};
              \node[main node, fill=white, label=above:{\tiny$-12$}, minimum size=8pt] (x2) at (0.7,0) {\scriptsize 2};
              \node[main node, fill=black, text=white, minimum size=6pt] (x3) at (1.4,0) {\scriptsize 3};
              \node[main node, fill=none, label=above:{\tiny$10$}, minimum size=9pt] (x3333) at (1.4,0) {};
              \node[main node, fill=none, dashed, shape=rectangle, minimum width=15pt, minimum height=12pt] (x33) at (1.4,0) {};
              \node[main node, fill=black, text=white, label=below:{\tiny$2$}, minimum size=8pt] (x4) at (0,-1) {\scriptsize 4};
              \node[main node, fill=black, text=white, minimum size=6pt] (x5) at (0.7,-1) {\scriptsize 5};
              \node[main node, fill=none, dashed, shape=rectangle, minimum width=52pt, minimum height=12pt] (x45) at (0.7,-1) {};
              \node[main node, fill=none, label=below:{\tiny$2$}, minimum size=9pt] (x5555) at (0.7,-1) {};
              \node[main node, fill=black, text=white, label=below:{\tiny$2$}, minimum size=8pt] (x6) at (1.4,-1) {\scriptsize 6};
              \node[main node, fill=none, dashed, shape=rectangle, minimum width=35pt, minimum height=12pt] (x45) at (0.35,0) {};
                \draw[->,red] (x1111) -- (x2);
                \draw (x3333) -- (x2);
                \draw (x1111) -- (x4);
                \draw[->,red] (x5555) -- (x4);
                \draw[->,red] (x5555) -- (x6);
                \draw (x6) -- (x3333);
            \end{tikzpicture}
            \caption{{ $3^{\textup{rd}}$ step}}
        \end{subfigure}
        \hspace{5mm}
        \begin{subfigure}{0.1\textwidth}
        \centering
            \begin{tikzpicture}[>=stealth',shorten >=1pt,auto,node distance=1.5cm,
            thick,main node/.style={circle,fill=black,draw,minimum size=5pt,inner sep=0pt}]
              \node[main node, fill=white, text=black, minimum size=6pt] (x1) at (0,0) {\scriptsize 1};
              \node[main node, fill=none, label=above:{\tiny$-2$}, minimum size=9pt] (x1111) at (0,0) {};
              \node[main node, fill=white, label=above:{\tiny$-2$}, minimum size=8pt] (x2) at (0.7,0) {\scriptsize 2};
              \node[main node, fill=white, text=black, minimum size=6pt] (x3) at (1.4,0) {\scriptsize 3};
              \node[main node, fill=none, label=above:{\tiny$-2$}, minimum size=9pt] (x3333) at (1.4,0) {};
              \node[main node, fill=black, text=white, label=below:{\tiny$2$}, minimum size=8pt] (x4) at (0,-1) {\scriptsize 4};
              \node[main node, fill=black, text=white, minimum size=6pt] (x5) at (0.7,-1) {\scriptsize 5};
              \node[main node, fill=none, dashed, shape=rectangle, minimum width=52pt, minimum height=12pt] (x45) at (0.7,-1) {};
              \node[main node, fill=none, label=below:{\tiny$2$}, minimum size=9pt] (x5555) at (0.7,-1) {};
              \node[main node, fill=black, text=white, label=below:{\tiny$2$}, minimum size=8pt] (x6) at (1.4,-1) {\scriptsize 6};
              \node[main node, fill=none, dashed, shape=rectangle, minimum width=52pt, minimum height=12pt] (x45) at (0.7,0) {};
                \draw[->,red] (x1111) -- (x2);
                \draw[->,red] (x3333) -- (x2);
                \draw (x1111) -- (x4);
                \draw[->,red] (x5555) -- (x4);
                \draw[->,red] (x5555) -- (x6);
                \draw (x6) -- (x3333);
            \end{tikzpicture}
            \caption{{ $4^{\textup{th}}$ step}}
        \end{subfigure}
        \hspace{5mm}
        \begin{subfigure}{0.1\textwidth}
        \centering
            \begin{tikzpicture}[>=stealth',shorten >=1pt,auto,node distance=1.5cm,
            thick,main node/.style={circle,fill=black,draw,minimum size=5pt,inner sep=0pt}]
              \node[main node, fill=white, text=black, minimum size=6pt] (x1) at (0,0) {\scriptsize 1};
              \node[main node, fill=none, label=above:{\tiny$0$}, minimum size=9pt] (x1111) at (0,0) {};
              \node[main node, fill=white, label=above:{\tiny$0$}, minimum size=8pt] (x2) at (0.7,0) {\scriptsize 2};
              \node[main node, fill=white, text=black, minimum size=6pt] (x3) at (1.4,0) {\scriptsize 3};
              \node[main node, fill=none, label=above:{\tiny$0$}, minimum size=9pt] (x3333) at (1.4,0) {};
              \node[main node, fill=white, text=black, label=below:{\tiny$0$}, minimum size=8pt] (x4) at (0,-1) {\scriptsize 4};
              \node[main node, fill=white, text=black, minimum size=6pt] (x5) at (0.7,-1) {\scriptsize 5};
              \node[main node, fill=none, label=below:{\tiny$0$}, minimum size=9pt] (x5555) at (0.7,-1) {};
              \node[main node, fill=white, text=black, label=below:{\tiny$0$}, minimum size=8pt] (x6) at (1.4,-1) {\scriptsize 6};
              \draw[dashed] (-0.2,0.2) -- (1.6,0.2);
              \draw[dashed] (-0.2,-1.2) -- (1.6,-1.2);
              \draw[dashed] (1.6,-1.2) -- (1.6,0.2);
              \draw[dashed] (-0.2,0.2) -- (-0.2,-0.25);
              \draw[dashed] (-0.2,-1.2) -- (-0.2,-0.8);
              \draw[dashed] (-0.2,-0.25) -- (1.2,-0.25);
              \draw[dashed] (-0.2,-0.8) -- (1.2,-0.8);
              \draw[dashed] (1.2,-0.25) -- (1.2,-0.8);
                \draw[->,red] (x1111) -- (x2);
                \draw[->,red] (x3333) -- (x2);
                \draw (x1111) -- (x4);
                \draw[->,red] (x5555) -- (x4);
                \draw[->,red] (x5555) -- (x6);
                \draw[->,red] (x6) -- (x3333);
            \end{tikzpicture}
            \caption{{ $5^{\textup{th}}$ step}}
        \end{subfigure}
    \caption{{\footnotesize Illustration of the sampling procedure. In b) and c), note that prioritization in \textsf{Sampler} avoids the flow blockage.}}
    \label{fig::feasibility}
\end{figure}
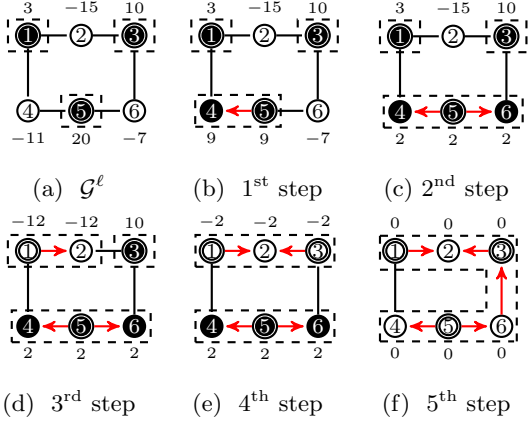
Figure~\ref{fig::feasibility} illustrates an example case where, after repeated execution of the \textsf{Sampler} function, \texttt{FORWARD} results in a feasible radial distribution in the shown subgraph $\mathcal{G}^\ell$. This is an example where we do not have capacity constraints on the edges.
{\small
\begin{algorithm}[t]
\footnotesize
\caption{\textsf{Sampler}}
\begin{algorithmic}[1]
\Require Concatenated dual graph $\bar{\mathcal{G}}^\ell$, existing polytrees $\mathbb{T}^\ell$, nodal value vector $\vect{p}^\ell$ 
\State $q\gets\emptyset$; $\bar{q}\gets\emptyset$; $\hat{q}\gets\emptyset$
\State $q\gets\textsf{Weight}(\bar{\mathcal{G}}^\ell, \vect{p}^\ell)$
\State Sort $q$ in decreasing order of weights
\While{$q$ is not empty}
    \State $((u,v,i,j), w_{i,j})\gets\textsf{Retrieve}(q)$ \Comment{See Appendix~\ref{Appedix::subroutines}}
    \If{$v$ is a pendent un-sampled super node in $\bar{\mathcal{G}}^\ell$}
        \State $e\gets(i\to j)$ 
        \State \textbf{Break}
    \ElsIf {$\vect{p}_i^\ell>0$} 
        \If {$\min\{\vect{p}_i^\ell,x_{i,j}\}+\min\{\vect{p}_j^\ell,\bar{\vect{p}}_v^\ell\}>0$}
        \State $\mathsf{Insert}((i\to j),\bar{q})$ 
        \Else
        \State $\mathsf{Insert}((i\to j),\hat{q})$ 
        \EndIf
        \Else
        \If {$\vect{p}_i^\ell+\min\{\vect{p}_j^\ell,\bar{x}_{i,i}\}>0$}
        \State $\mathsf{Insert}((j\to i),\bar{q})$ 
        \Else 
        \State $\mathsf{Insert}((j\to i),\hat{q})$ 
        \EndIf
    \EndIf
\EndWhile
\State $e\gets \textsf{Retrieve}\left([\bar{q},\hat{q}]\right)$
\State \Return $e$
\end{algorithmic}
\label{alg::sampler}
\end{algorithm}
}
{\small 
\begin{algorithm}[t]
{\footnotesize
\caption{\textsf{Rewire}}
\begin{algorithmic}[1]
\Require Sub-graph partition $\bar{\mathcal{G}}^\ell$, flow vector $\vect{p}^\ell$ and $\mathcal{S}^\ell$
\State Find $v \leftarrow \{v\in\mathcal{G}^\ell|\vect{p}^\ell_v<0\}$
\State Identify $\mathcal{C}(\mathcal{R},v)$ from $\mathcal{N}(v)\in\mathcal{S}^\ell$
\State Define $\mathcal{C}(\mathcal{R}^+)$ by finding all nodes with surplus flow in $\mathcal{S}^\ell$
\State Define $\mathcal{C}(\bar{\mathcal{R}}) \leftarrow \{\mathcal{R}\in\mathcal{S}^\ell | \vect{p}^\ell_{u^\prime}<\bar{x}(u^\prime,v^\prime),\quad\forall(u^\prime,v^\prime)\in\mathcal{R}\}$
\For{$\mathcal{R}_{j}^{i}\in\mathcal{C}(\mathcal{R},v)\setminus\mathcal{C}(\bar{\mathcal{R}})$}
    \State Find $\{\mathcal{R}_{j}^{i^\prime}\in\mathcal{C}(\mathcal{R}^+)|\exists(u^\prime,v^\prime)\in\mathcal{G}^\ell\}$ where $u^\prime$ is leaf of $\mathcal{R}_{j}^{i^\prime}$ and $v^\prime\in\mathcal{T}^\ell_i$.
    \State Find $\hat{e}_{u^\prime,v^\prime}$ where $u^\prime\in\mathcal{R}_{j}^{i}$ and $v^\prime\notin\mathcal{R}_{j}^{i}$ and $v^\prime\in\mathcal{T}^\ell_{i}$
    \If{$\exists \mathcal{R}_{j}^{i^\prime}$}
    \State Delete $\hat{e}_{u^\prime,v^\prime}$
    \State Define $\mathcal{R}_{j}^{i}\leftarrow\mathcal{R}_{j}^{i}\cup\{v\}$, $\mathcal{R}_{j}^{i}\leftarrow\mathcal{R}_{j}^{i}\setminus\{v^\prime\}$, $\mathcal{R}_{j}^{i^\prime}\leftarrow\mathcal{R}_{j}^{i}\cup\{v^\prime\}$
    \State $break$
    \EndIf
\EndFor
\State \Return $\mathcal{S}^\ell$
\end{algorithmic}
\label{alg::Rewire}
}
\end{algorithm}
}
\subsubsection{\textsf{Rewire}}\label{sec::rewire}
\vspace{-0.1in}
When capacity constraints prevent \textsf{Sampler} from supplying all nodes, \textsf{Rewire} (Algorithm~\ref{alg::Rewire}) corrects the infeasible solution through strategic edge swapping. The key insight is that unsupplied nodes can be reached through alternative paths if we reroute flow from oversupplied~regions.

We define a \emph{source-to-leaf path} $\mathcal{R}_{k}^{i}$ as a path from a source node to a leaf node within polytree $\mathcal{T}^\ell_i$. Note that $\mathcal{T}^\ell_i = \bigcup_{k=1}^{r}\mathcal{R}_{k}^{i}$, where $r$ is the total number of the leaf nodes in polytree $\mathcal{T}^\ell_i$. The algorithm uses three key sets: $\mathcal{C}(\mathcal{R},v)$ contains source-to-leaf paths that have edges connecting to unsupplied node $v$ in the original graph $\mathcal{G}^\ell$; $\mathcal{C}(\mathcal{R}^+)$ contains \textit{surplus source-to-leaf paths} from the already-sampled polytrees where leaf nodes have remaining supply ($\vect{p}^\ell > 0$); and $\mathcal{C}(\bar{\mathcal{R}})$ contains \textit{saturated source-to-leaf paths} from the sampled polytrees where leaf nodes are not fully supplied due to capacity limits along the~path.

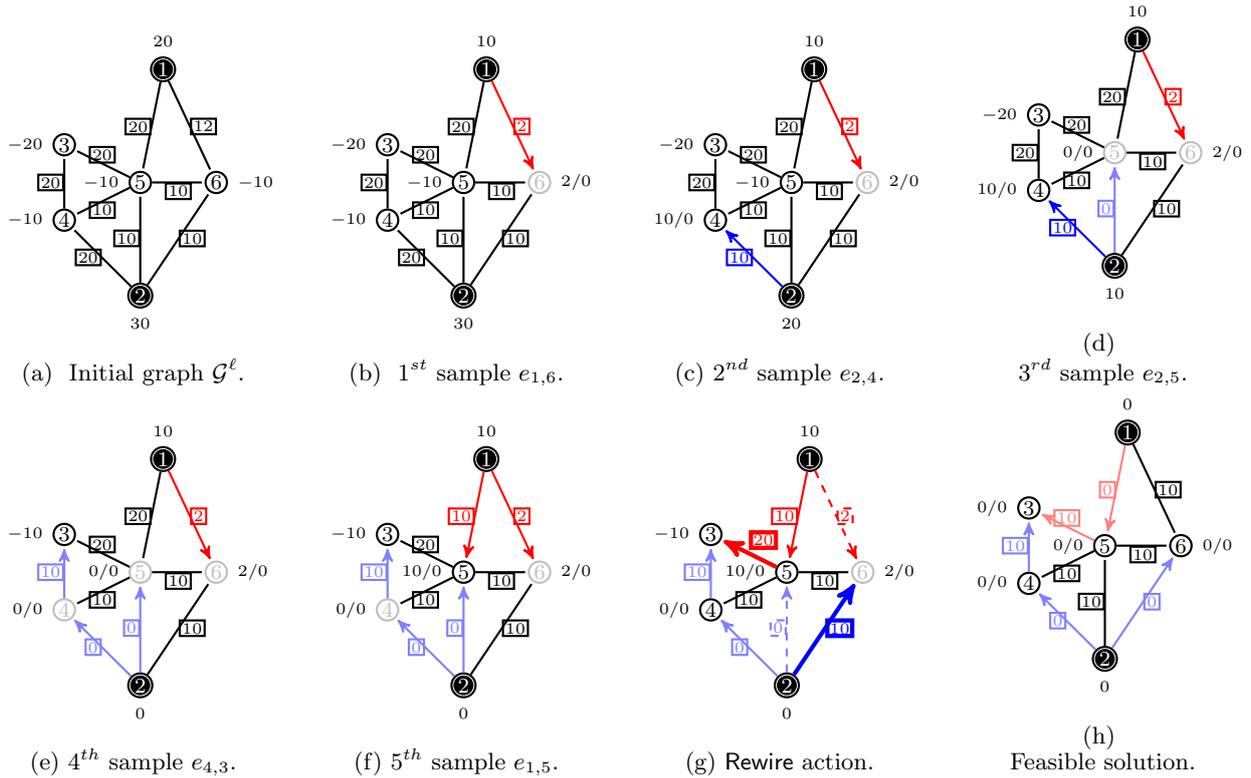
\begin{figure*}[t]
    \centering
    \begin{subfigure}{0.2\textwidth}
            \centering
            \begin{tikzpicture}[>=stealth',shorten      >=1pt,auto,node distance=1.5cm,
            thick,main node/.style={circle,fill=black,draw,minimum size=5pt,inner sep=0pt}]
              \node[main node, fill=black, text=white, minimum size=6pt] (x1) at (1.3,1.5) {\scriptsize 1};
              \node[main node, fill=none, label=above:{\tiny$20$}, minimum size=9pt] (x1111) at (1.3,1.5) {};
              \node[main node, fill=black, text=white, minimum size=6pt] (x2) at (1,-1.5) {\scriptsize 2};
              \node[main node, fill=none, label=below:{\tiny$30$}, minimum size=9pt] (x2222) at (1,-1.5) {};
              \node[main node, fill=white, label=left:{\tiny$-20$}, minimum size=8pt] (x3) at (0,0.5) {\scriptsize 3};
              \node[main node, fill=white, label=left:{\tiny$-10$}, minimum size=8pt] (x4) at (0,-0.5) {\scriptsize 4};
              \node[main node, fill=white, label=left:{\tiny$-10$}, minimum size=8pt] (x5) at (1,0) {\scriptsize 5};
              \node[main node, fill=white, label=right:{\tiny$-10$}, minimum size=8pt] (x6) at (2,0) {\scriptsize 6};
                \draw  (x1111) -- node[left,draw,inner sep=1pt] {\tiny $20$} (x5);
                \draw  (x1111) -- node[right,draw,inner sep=1pt] {\tiny $12$} (x6);
                \draw  (x5) -- node[above,draw,inner sep=1pt] {\tiny $20$} (x3);
                \draw  (x5) -- node[below,draw,inner sep=1pt] {\tiny $10$} (x6);
                \draw  (x5) -- node[below,draw,inner sep=1pt] {\tiny $10$} (x4);
                \draw  (x4) -- node[left,draw,inner sep=1pt] {\tiny $20$} (x3);
                \draw  (x2222) -- node[left,draw,inner sep=1pt] {\tiny $10$} (x5);
                \draw  (x2222) -- node[right,draw,inner sep=1pt] {\tiny $10$} (x6);
                \draw  (x2222) -- node[left,draw,inner sep=1pt] {\tiny $20$} (x4);
            \end{tikzpicture}
            \caption{{ Initial graph $\mathcal{G}^\ell$.}}
        \end{subfigure}
        \hspace{5mm}
    \begin{subfigure}{0.2\textwidth}
            \centering
            \begin{tikzpicture}[>=stealth',shorten      >=1pt,auto,node distance=1.5cm,
            thick,main node/.style={circle,fill=black,draw,minimum size=5pt,inner sep=0pt}]
              \node[main node, fill=black, text=white, minimum size=6pt] (x1) at (1.3,1.5) {\scriptsize 1};
              \node[main node, fill=none, label=above:{\tiny$10$}, minimum size=9pt] (x1111) at (1.3,1.5) {};
              \node[main node, fill=black, text=white, minimum size=6pt] (x2) at (1,-1.5) {\scriptsize 2};
              \node[main node, fill=none, label=below:{\tiny$30$}, minimum size=9pt] (x2222) at (1,-1.5) {};
              \node[main node, fill=white, label=left:{\tiny$-20$}, minimum size=8pt] (x3) at (0,0.5) {\scriptsize 3};
              \node[main node, fill=white, label=left:{\tiny$-10$}, minimum size=8pt] (x4) at (0,-0.5) {\scriptsize 4};
              \node[main node, fill=white, label=left:{\tiny$-10$}, minimum size=8pt] (x5) at (1,0) {\scriptsize 5};
              \node[main node, fill=white, draw=gray!50, label=right:{\tiny$2/0$}, minimum size=8pt] (x6) at (2,0) {\scriptsize {\color{gray!50}6}};
                \draw  (x1111) -- node[left,draw,inner sep=1pt] {\tiny $20$} (x5);
                \draw[->,red]  (x1111) -- node[right,draw,inner sep=1pt] {\tiny $2$} (x6);
                \draw  (x5) -- node[above,draw,inner sep=1pt] {\tiny $20$} (x3);
                \draw  (x5) -- node[below,draw,inner sep=1pt] {\tiny $10$} (x6);
                \draw  (x5) -- node[below,draw,inner sep=1pt] {\tiny $10$} (x4);
                \draw  (x4) -- node[left,draw,inner sep=1pt] {\tiny $20$} (x3);
                \draw  (x2222) -- node[left,draw,inner sep=1pt] {\tiny $10$} (x5);
                \draw  (x2222) -- node[right,draw,inner sep=1pt] {\tiny $10$} (x6);
                \draw  (x2222) -- node[left,draw,inner sep=1pt] {\tiny $20$} (x4);
            \end{tikzpicture}
            \caption{{ $1^{st}$ sample $e_{1,6}$.}}
        \end{subfigure}
        \hspace{5mm}
    \begin{subfigure}{0.2\textwidth}
            \centering
            \begin{tikzpicture}[>=stealth',shorten      >=1pt,auto,node distance=1.5cm,
            thick,main node/.style={circle,fill=black,draw,minimum size=5pt,inner sep=0pt}]
              \node[main node, fill=black, text=white, minimum size=6pt] (x1) at (1.3,1.5) {\scriptsize 1};
              \node[main node, fill=none, label=above:{\tiny$10$}, minimum size=9pt] (x1111) at (1.3,1.5) {};
              \node[main node, fill=black, text=white, minimum size=6pt] (x2) at (1,-1.5) {\scriptsize 2};
              \node[main node, fill=none, label=below:{\tiny$20$}, minimum size=9pt] (x2222) at (1,-1.5) {};
              \node[main node, fill=white, label=left:{\tiny$-20$}, minimum size=8pt] (x3) at (0,0.5) {\scriptsize 3};
              \node[main node, fill=white, label=left:{\tiny$10/0$}, minimum size=8pt] (x4) at (0,-0.5) {\scriptsize 4};
              \node[main node, fill=white, label=left:{\tiny$-10$}, minimum size=8pt] (x5) at (1,0) {\scriptsize 5};
              \node[main node, draw=gray!50, fill=white, label=right:{\tiny$2/0$}, minimum size=8pt] (x6) at (2,0) {\scriptsize \textcolor{gray!50}{6}};
                \draw  (x1111) -- node[left,draw,inner sep=1pt] {\tiny $20$} (x5);
                \draw[->,red]  (x1111) -- node[right,draw,inner sep=1pt] {\tiny $2$} (x6);
                \draw  (x5) -- node[above,draw,inner sep=1pt] {\tiny $20$} (x3);
                \draw  (x5) -- node[below,draw,inner sep=1pt] {\tiny $10$} (x6);
                \draw  (x5) -- node[below,draw,inner sep=1pt] {\tiny $10$} (x4);
                \draw  (x4) -- node[left,draw,inner sep=1pt] {\tiny $20$} (x3);
                \draw  (x2222) -- node[left,draw,inner sep=1pt] {\tiny $10$} (x5);
                \draw  (x2222) -- node[right,draw,inner sep=1pt] {\tiny $10$} (x6);
                \draw[->,blue]  (x2222) -- node[left,draw,inner sep=1pt] {\tiny $10$} (x4);
            \end{tikzpicture}
            \caption{{$2^{nd}$ sample $e_{2,4}$.}}
        \end{subfigure}
        \hspace{5mm}
    \begin{subfigure}{0.2\textwidth}
            \centering
            \begin{tikzpicture}[>=stealth',shorten      >=1pt,auto,node distance=1.5cm,
            thick,main node/.style={circle,fill=black,draw,minimum size=5pt,inner sep=0pt}]
              \node[main node, fill=black, text=white, minimum size=6pt] (x1) at (1.3,1.5) {\scriptsize 1};
              \node[main node, fill=none, label=above:{\tiny$10$}, minimum size=9pt] (x1111) at (1.3,1.5) {};
              \node[main node, fill=black, text=white, minimum size=6pt] (x2) at (1,-1.5) {\scriptsize 2};
              \node[main node, fill=none, label=below:{\tiny$10$}, minimum size=9pt] (x2222) at (1,-1.5) {};
              \node[main node, fill=white, label=left:{\tiny$-20$}, minimum size=8pt] (x3) at (0,0.5) {\scriptsize 3};
              \node[main node, fill=white, label=left:{\tiny$10/0$}, minimum size=8pt] (x4) at (0,-0.5) {\scriptsize 4};
              \node[main node, fill=white, draw=gray!50,label=left:{\tiny$0/0$}, minimum size=8pt] (x5) at (1,0) {\scriptsize {\color{gray!50}5}};
              \node[main node, fill=white, draw=gray!50,label=right:{\tiny$2/0$}, minimum size=8pt] (x6) at (2,0) {\scriptsize {\color{gray!50}6}};
                \draw  (x1111) -- node[left,draw,inner sep=1pt] {\tiny $20$} (x5);
                \draw[->,red]  (x1111) -- node[right,draw,inner sep=1pt] {\tiny $2$} (x6);
                \draw  (x5) -- node[above,draw,inner sep=1pt] {\tiny $20$} (x3);
                \draw  (x5) -- node[below,draw,inner sep=1pt] {\tiny $10$} (x6);
                \draw  (x5) -- node[below,draw,inner sep=1pt] {\tiny $10$} (x4);
                \draw  (x4) -- node[left,draw,inner sep=1pt] {\tiny $20$} (x3);
                \draw[->,blue!50!white]  (x2222) -- node[left,draw,inner sep=1pt] {\tiny $0$} (x5);
                \draw  (x2222) -- node[right,draw,inner sep=1pt] {\tiny $10$} (x6);
                \draw[->,blue]  (x2222) -- node[left,draw,inner sep=1pt] {\tiny $10$} (x4);
            \end{tikzpicture}
            \caption{{\small $3^{rd}$ sample $e_{2,5}$.}}
        \end{subfigure}
        \hspace{5mm}
    \begin{subfigure}{0.2\textwidth}
            \centering
            \begin{tikzpicture}[>=stealth',shorten      >=1pt,auto,node distance=1.5cm,
            thick,main node/.style={circle,fill=black,draw,minimum size=5pt,inner sep=0pt}]
              \node[main node, fill=black, text=white, minimum size=6pt] (x1) at (1.3,1.5) {\scriptsize 1};
              \node[main node, fill=none, label=above:{\tiny$10$}, minimum size=9pt] (x1111) at (1.3,1.5) {};
              \node[main node, fill=black, text=white, minimum size=6pt] (x2) at (1,-1.5) {\scriptsize 2};
              \node[main node, fill=none, label=below:{\tiny$0$}, minimum size=9pt] (x2222) at (1,-1.5) {};
              \node[main node, fill=white, label=left:{\tiny$-10$}, minimum size=8pt] (x3) at (0,0.5) {\scriptsize 3};
              \node[main node, draw=gray!50, fill=white, label=left:{\tiny$0/0$}, minimum size=8pt] (x4) at (0,-0.5) {\scriptsize \color{gray!50} 4};
              \node[main node, draw=gray!50, fill=white, label=left:{\tiny$0/0$}, minimum size=8pt] (x5) at (1,0) {\scriptsize \color{gray!50} 5};
              \node[main node, draw=gray!50, fill=white, label=right:{\tiny$2/0$}, minimum size=8pt] (x6) at (2,0) {\scriptsize \color{gray!50} 6};
                \draw  (x1111) -- node[left,draw,inner sep=1pt] {\tiny $20$} (x5);
                \draw[->,red]  (x1111) -- node[right,draw,inner sep=1pt] {\tiny $2$} (x6);
                \draw  (x5) -- node[above,draw,inner sep=1pt] {\tiny $20$} (x3);
                \draw  (x5) -- node[below,draw,inner sep=1pt] {\tiny $10$} (x6);
                \draw  (x5) -- node[below,draw,inner sep=1pt] {\tiny $10$} (x4);
                \draw[->,blue!50!white]  (x4) -- node[left,draw,inner sep=1pt] {\tiny $10$} (x3);
                \draw[->,blue!50!white]  (x2222) -- node[left,draw,inner sep=1pt] {\tiny $0$} (x5);
                \draw  (x2222) -- node[right,draw,inner sep=1pt] {\tiny $10$} (x6);
                \draw[->,blue!50!white]  (x2222) -- node[left,draw,inner sep=1pt] {\tiny $0$} (x4);
            \end{tikzpicture}
            \caption{{$4^{th}$ sample $e_{4,3}$.}}
        \end{subfigure}
        \hspace{5mm}
    \begin{subfigure}{0.2\textwidth}
            \centering
            \begin{tikzpicture}[>=stealth',shorten      >=1pt,auto,node distance=1.5cm,
            thick,main node/.style={circle,fill=black,draw,minimum size=5pt,inner sep=0pt}]
              \node[main node, fill=black, text=white, minimum size=6pt] (x1) at (1.3,1.5) {\scriptsize 1};
              \node[main node, fill=none, label=above:{\tiny$10$}, minimum size=9pt] (x1111) at (1.3,1.5) {};
              \node[main node, fill=black, text=white, minimum size=6pt] (x2) at (1,-1.5) {\scriptsize 2};
              \node[main node, fill=none, label=below:{\tiny$0$}, minimum size=9pt] (x2222) at (1,-1.5) {};
              \node[main node, fill=white, label=left:{\tiny$-10$}, minimum size=8pt] (x3) at (0,0.5) {\scriptsize 3};
              \node[main node, color=gray!50, fill=white, label=left:{\tiny  $0/0$}, minimum size=8pt] (x4) at (0,-0.5) {\scriptsize \color{gray!50} 4};
              \node[main node, fill=white, label=left:{\tiny$10/0$}, minimum size=8pt] (x5) at (1,0) {\scriptsize 5};
              \node[main node, fill=white, draw=gray!50,label=right:{\tiny$2/0$}, minimum size=8pt] (x6) at (2,0) {\scriptsize {\color{gray!50}6}};
                \draw[->,red]  (x1111) -- node[left,draw,inner sep=1pt] {\tiny $10$} (x5);
                \draw[->,red]  (x1111) -- node[right,draw,inner sep=1pt] {\tiny $2$} (x6);
                \draw  (x5) -- node[above,draw,inner sep=1pt] {\tiny $20$} (x3);
                \draw  (x5) -- node[below,draw,inner sep=1pt] {\tiny $10$} (x6);
                \draw  (x5) -- node[below,draw,inner sep=1pt] {\tiny $10$} (x4);
                \draw[->,blue!50!white]  (x4) -- node[left,draw,inner sep=1pt] {\tiny $10$} (x3);
                \draw[->,blue!50!white]  (x2222) -- node[left,draw,inner sep=1pt] {\tiny $0$} (x5);
                \draw  (x2222) -- node[right,draw,inner sep=1pt] {\tiny $10$} (x6);
                \draw[->,blue!50!white]  (x2222) -- node[left,draw,inner sep=1pt] {\tiny $0$} (x4);
            \end{tikzpicture}
            \caption{{$5^{th}$ sample $e_{1,5}$.}}
        \end{subfigure}
        \hspace{5mm}
        \begin{subfigure}{0.2\textwidth}
            \centering
            \begin{tikzpicture}[>=stealth',shorten      >=1pt,auto,node distance=1.5cm,
            thick,main node/.style={circle,fill=black,draw,minimum size=5pt,inner sep=0pt}]
              \node[main node, fill=black, text=white, minimum size=6pt] (x1) at (1.3,1.5) {\scriptsize 1};
              \node[main node, fill=none, label=above:{\tiny$10$}, minimum size=9pt] (x1111) at (1.3,1.5) {};
              \node[main node, fill=black, text=white, minimum size=6pt] (x2) at (1,-1.5) {\scriptsize 2};
              \node[main node, fill=none, label=below:{\tiny$0$}, minimum size=9pt] (x2222) at (1,-1.5) {};
              \node[main node, fill=white, label=left:{\tiny$-10$}, minimum size=8pt] (x3) at (0,0.5) {\scriptsize 3};
              \node[main node, fill=white, label=left:{\tiny$0/0$}, minimum size=8pt] (x4) at (0,-0.5) {\scriptsize 4};
              \node[main node, fill=white, label=left:{\tiny$10/0$}, minimum size=8pt] (x5) at (1,0) {\scriptsize 5};
              \node[main node, draw=gray!50, fill=white, label=right:{\tiny$2/0$}, minimum size=8pt] (x6) at (2,0) {\scriptsize \color{gray!50} 6};
                \draw[->,red]  (x1111) -- node[left,draw,inner sep=1pt] {\tiny $10$} (x5);
                \draw[->,red, dashed]  (x1111) -- node[right,draw,inner sep=1pt] {\tiny $2$} (x6);
                \draw[->,red,ultra thick]  (x5) -- node[right,yshift=5pt,draw,inner sep=1pt] {\tiny $20$} (x3);
                \draw  (x5) -- node[below,draw,inner sep=1pt] {\tiny $10$} (x6);
                \draw  (x5) -- node[below,draw,inner sep=1pt] {\tiny $10$} (x4);
                \draw[->,blue!50]  (x4) -- node[left,draw,inner sep=1pt] {\tiny $10$} (x3);
                \draw[->,blue!50!white,dashed]  (x2222) -- node[left,draw,inner sep=1pt] {\tiny $0$} (x5);
                \draw[->,blue,ultra thick]  (x2222) -- node[right,draw,inner sep=1pt] {\tiny $10$} (x6);
                \draw[->,blue!50!white]  (x2222) -- node[left,draw,inner sep=1pt] {\tiny $0$} (x4);
            \end{tikzpicture}
            \caption{{\textsf{Rewire} action.}}
        \end{subfigure}
        \hspace{5mm}
        \begin{subfigure}{0.2\textwidth}
            \centering
            \begin{tikzpicture}[>=stealth',shorten      >=1pt,auto,node distance=1.5cm,
            thick,main node/.style={circle,fill=black,draw,minimum size=5pt,inner sep=0pt}]
              \node[main node, fill=black, text=white, minimum size=6pt] (x1) at (1.3,1.5) {\scriptsize 1};
              \node[main node, fill=none, label=above:{\tiny$0$}, minimum size=9pt] (x1111) at (1.3,1.5) {};
              \node[main node, fill=black, text=white, minimum size=6pt] (x2) at (1,-1.5) {\scriptsize 2};
              \node[main node, fill=none, label=below:{\tiny$0$}, minimum size=9pt] (x2222) at (1,-1.5) {};
              \node[main node, fill=white, label=left:{\tiny$0/0$}, minimum size=8pt] (x3) at (0,0.5) {\scriptsize 3};
              \node[main node, fill=white, label=left:{\tiny$0/0$}, minimum size=8pt] (x4) at (0,-0.5) {\scriptsize 4};
              \node[main node, fill=white, label=left:{\tiny$0/0$}, minimum size=8pt] (x5) at (1,0) {\scriptsize 5};
              \node[main node, fill=white, label=right:{\tiny$0/0$}, minimum size=8pt] (x6) at (2,0) {\scriptsize 6};
                \draw[->,red!50!white]  (x1111) -- node[left,draw,inner sep=1pt] {\tiny $0$} (x5);
                \draw  (x1111) -- node[right,draw,inner sep=1pt] {\tiny $10$} (x6);
                \draw[->,red!50!white]  (x5) -- node[above,draw,inner sep=1pt] {\tiny $10$} (x3);
                \draw  (x5) -- node[below,draw,inner sep=1pt] {\tiny $10$} (x6);
                \draw  (x5) -- node[below,draw,inner sep=1pt] {\tiny $10$} (x4);
                \draw[->,blue!50!white]  (x4) -- node[left,draw,inner sep=1pt] {\tiny $10$} (x3);
                \draw  (x2222) -- node[left,draw,inner sep=1pt] {\tiny $10$} (x5);
                \draw[->,blue!50!white]  (x2222) -- node[right,draw,inner sep=1pt] {\tiny $0$} (x6);
                \draw[->,blue!50!white]  (x2222) -- node[left,draw,inner sep=1pt] {\tiny $0$} (x4);
            \end{tikzpicture}
            \caption{{\small Feasible solution.}}
        \end{subfigure}
    \caption{{\footnotesize An example demonstrating the application of the \textsf{Rewire} function to resolve an infeasible solution generated by \textsf{Sampler} alone. Input and output values at each node in the original graph $\mathcal{G}^\ell$ are shown as numbers next to the nodes. Edge capacity values are displayed in boxed weights. After sampling, the notation $x/y$ next to each sink node indicates: $x$ represents the excess supply available at that node, and $y$ represents the demand at the node. When an edge is sampled to supply downstream nodes, the capacity values are adjusted based on the remaining flow that can pass through the edge, to be the excess supply at the terminal node of the sampled edge.  
    }}
    \label{fig::unfeasible_capacity}
\end{figure*}
\textsf{Rewire} operates through a systematic swap mechanism to redistribute flow from oversupplied regions to unsupplied nodes. The algorithm works in three main phases:

\textbf{Phase 1: Identify rewiring opportunities.} The algorithm first locates unsupplied nodes $v$ where $\vect{p}^\ell_v < 0$ (line 1) and identifies source-to-leaf paths that could potentially supply these nodes through non-sampled edges in the original graph, specifically $\mathcal{C}(\mathcal{R},v)\setminus\mathcal{C}(\bar{\mathcal{R}})$ (lines 2-4). These represent source-to-leaf paths with sufficient capacity to reach the unsupplied nodes if properly connected. Simultaneously, the algorithm identifies surplus source-to-leaf paths $\mathcal{C}(\mathcal{R}^+)$ whose leaf nodes have remaining supply that can be redistributed.

\textbf{Phase 2: Execute strategic swaps.} For each viable source-to-leaf path in $\mathcal{C}(\mathcal{R},v)\setminus\mathcal{C}(\bar{\mathcal{R}})$, the algorithm seeks surplus source-to-leaf paths $\mathcal{C}(\mathcal{R}^+)$ whose leaf nodes can donate flow (lines 5-6). When such a surplus source-to-leaf path is found, \textsf{Rewire} performs a swap operation: it removes a blocking edge $\hat{e}_{u',v'}$ from the current polytree structure (line 9) and establishes new connections that redirect the surplus flow to supply the previously unsupplied node $v$ (line 10). This swap effectively transfers the surplus leaf node from one polytree to another, creating a new supply path to the unsupplied node.

\textbf{Phase 3: Maintain radiality.} The edge deletion and reconnection process ensures that the resulting network maintains its radial (tree-like) structure by avoiding cycles while establishing new supply paths to previously unsupplied nodes (lines 9-11). The algorithm manages these topological changes to preserve the polytree properties essential for the overall solution framework.

The key insight is that \textsf{Rewire} leverages the flexibility of polytree structures to redistribute existing supply capacity rather than requiring additional generation resources. By strategically swapping connections between surplus and deficit regions, the algorithm transforms an infeasible solution into a feasible one, as illustrated in Fig.~\ref{fig::unfeasible_capacity}. This approach is particularly effective because it maintains the computational efficiency of working within the sampled polytree structure while resolving capacity-induced infeasibilities through localized network modifications.

\section{Feasibility Analysis of FORWARD Generated Radial Distribution Networks}
\label{sec::feasibility}
\vspace{-0.1in}
This section demonstrates that the \texttt{FORWARD} algorithm is guaranteed to generate feasible radial distribution configurations for optimization problem~\eqref{eqn::problem1}. By Theorem~\ref{thm::partition}, we established that feasibility of the original problem is equivalent to finding feasible solutions in each balanced component $\mathcal{G}^\ell$. Therefore, our analysis focuses on proving that \texttt{FORWARD} generates feasible radial configurations within each subgraph $\mathcal{G}^\ell$.

\subsection{Feasibility Without Capacity Constraints}
\vspace{-0.1in}
We begin by analyzing the case where capacity constraints are not binding, establishing the fundamental feasibility properties of \texttt{FORWARD}.

\begin{lem}\longthmtitle{Edge Direction Preservation in Sampled Polytrees}
    Given an edge $e$ sampled by \textsf{Sampler} to merge two distinct polytrees, it will not reverse the directionality of previously sampled edges.
    \label{lem::direction_preservation}
\end{lem}
\begin{pf}
    Assume for contradiction that merging polytrees $\mathcal{T}^\ell_i$ (with aggregate nodal value $\vect{p}^\ell_i > 0$) and $\mathcal{T}^\ell_j$ (with $\vect{p}^\ell_j < 0$) via edge $e$ causes a previously established flow direction within $\mathcal{T}^\ell_j$ to reverse. For flow reversal to occur within $\mathcal{T}^\ell_j$, there must have been an internal flow path directed towards a sink node $u$, which upon connection via $e$ is forced to redirect away from $u$. This would only happen if the merged polytree has net positive flow that redirects existing internal flows of $\mathcal{T}^\ell_j$.  However, the \textsf{Sampler}'s priority condition requires $\vect{p}^\ell_i + \vect{p}^\ell_j \geq 0$ for edge selection. This condition ensures that the merged polytree maintains sufficient supply to meet existing demands without forcing flow reversals. The quasi-bipartite property of $\bar{\mathcal{G}}^\ell$ ensures edges connect supply-side components to demand-side components, with flow naturally moving from positive to negative nodal values. Within the tree structure of polytrees, flow directions are uniquely determined by source-sink relationships once demands are fixed (as established in Lemma~\ref{lem::uniqueness}). The \textsf{Sampler}'s aggregate nodal value check and the tree-like nature of polytrees ensure that previously directed edges maintain their flow direction when polytrees are merged.\boxend
\end{pf}

\begin{lem}\longthmtitle{Properties of \texttt{FORWARD}'s Output}
    After \texttt{FORWARD} execution on each subgraph $\mathcal{G}^\ell(\mathcal{V}^\ell,\mathcal{E}^\ell)$, the resulting graph $\mathcal{G}(\mathcal{V}(\mathcal{S}^\ell),\mathcal{S}^\ell)$ is a radial configuration covering all nodes ($\mathcal{V}(\mathcal{S}^\ell)=\mathcal{V}^\ell$). If the initial input-output flow is balanced in $\mathcal{G}^\ell$, then every node's demand is met with $\vect{p}_i^\ell = \mathbf{0}$ for all $i \in \mathcal{V}^\ell$.
    \label{lem::forward_properties}
\end{lem}
\begin{pf}
The proof demonstrates that \texttt{FORWARD} consistently generates a radial configuration, ensures complete node coverage, and guarantees flow balance when the initial subgraph is balanced.

\emph{Radial Configuration (Polyforest Property):}
The process initiates with each source node forming a trivial polytree. The \texttt{FORWARD} algorithm's core loop iteratively adds edges via \textsf{Sampler}. The \textsf{Sampler} operates on the dual graph $\bar{\mathcal{G}}^\ell$ constructed by \textsf{Net-Concad}, where existing polytrees are condensed into super sampled nodes and un-sampled components into super un-sampled nodes. The \textsf{Sampler}'s edge selection ensures that each chosen edge $(s,t)$ either (a) connects an existing polytree to an un-sampled component or (b) merges two distinct existing polytrees. In both cases, a single edge is added between previously disconnected components, strictly preventing cycle formation and ensuring $\mathcal{G}(\mathcal{V}(\mathcal{S}^{\ell}),\mathcal{S}^{\ell})$ remains a polyforest.

\emph{Complete Node Coverage and Flow Balance:}
We prove that the generated polyforest completely covers the graph ($\mathcal{V}(\mathcal{S}^\ell)=\mathcal{V}^\ell$) and delivers inputs from sources to meet all demands within a finite number of steps.

Upon \textsf{Islander} application, each subgraph $\mathcal{G}^\ell$ has one or more nodes with positive nodal value, collectively balancing the total demand across $\mathcal{G}^\ell$. Since each $\mathcal{G}^\ell$ has minimum degree 2, the dual graph $\bar{\mathcal{G}}^\ell$ constructed by \textsf{Net-Concad} is consistently connected. This connectivity ensures every un-sampled super node in $\bar{\mathcal{G}}^\ell$ connects to at least one super sampled node.

At each step, as long as not all nodes are covered, there exists at least one super sampled node with positive nodal value (since total flow in $\mathcal{G}^\ell$ is balanced and unmet sink demands require corresponding excess supply in sampled polytrees). The \textsf{Sampler} prioritizes edges based on its priority list, systematically extending polytrees by connecting them to un-sampled nodes or merging distinct polytrees.

The \textsf{Weight} function assigns positive weight $w_{i,j}>0$ only to edges where precisely one connected super-node has positive nodal value ($\textsf{XOR}(\vect{p}_i^\ell>0, \vect{p}_j^\ell>0)=1$), guiding \textsf{Sampler}'s choices. By prioritizing connections that can be fully supplied ($\vect{p}^\ell_i + \vect{p}^\ell_j>0$), this process effectively distributes positive nodal values from sources towards un-sampled super nodes until flow is balanced.

This iterative process terminates because $|\mathcal{V}^\ell|$ is finite. Each step either adds a new node to the covered set or merges existing polytrees, ensuring monotonic progress towards full coverage. The process stops when all sampled polytrees have zero nodal value, which occurs only when all nodes are covered ($\mathcal{V}(\mathcal{S}^\ell)=\mathcal{V}^\ell$) and all demands are met, resulting in $\vect{p}_i^\ell = \mathbf{0}$ for all $i \in \mathcal{V}^\ell$.\boxend
\end{pf}

\begin{thm}\longthmtitle{\texttt{FORWARD} Generates Feasible Radial Configurations (Uncapacitated Case)}
    The \texttt{FORWARD} algorithm (Algorithm 5) without capacity constraints is guaranteed to construct feasible radial configurations for optimization problem \eqref{eqn::problem1}.
    \label{thm::feasibility_uncapacitated}
\end{thm}
\begin{pf}
A feasible configuration for problem \eqref{eqn::problem1} requires radial structure (polyforest), complete node coverage, and balanced flow distribution. The proof analyzes the sequential operations of \texttt{FORWARD} and properties of its~outputs with regards to these attributes.

\begin{enumerate}
    \item \textbf{\textsf{Pre-Processor} Operation:}
    The \textsf{Pre-Processor} separates radial subgraphs from the 2-core subgraph ($\mathcal{G}_P$) at articulation points, balancing nodal values to ensure input-output balance within each component. By Lemma~\ref{lem::uniqueness}, the polyforest created over these separated radial configurations is unique and meets all demands within their respective components.

    \item \textbf{\textsf{Islander} Operation:}
    The \textsf{Islander} function partitions the 2-core subgraph $\mathcal{G}_P$ at articulation points with positive nodal value, adjusting nodal values to ensure each newly created subgraph $\mathcal{G}^\ell$ is input-output balanced.

    \item \textbf{\texttt{FORWARD} Operation on each $\mathcal{G}^\ell$:}
    For each balanced subgraph $\mathcal{G}^\ell$, Lemma~\ref{lem::forward_properties} proves that \texttt{FORWARD} ensures: (a) radial configuration construction, (b) complete node coverage, and (c) flow balance with all demands met.

    \item \textbf{Overall Feasibility:}
    Each component processed by \textsf{Pre-Processor} and each subgraph $\mathcal{G}^\ell$ processed by \texttt{FORWARD} yields a balanced radial configuration. These configurations are joined at their respective articulation points, resulting in a complete radial and flow-balanced configuration for the entire $\mathcal{G}_D$ that satisfies all feasibility criteria for optimization~problem~\eqref{eqn::problem1}.\boxend
\end{enumerate}
\end{pf}

\subsection{Feasibility With Capacity Constraints}
\vspace{-0.1in}
When capacity constraints are present, \texttt{FORWARD} may initially produce infeasible solutions where some nodes cannot be fully supplied due to capacity bottlenecks along their source-to-leaf paths. As described earlier, the \textsf{Rewire} function addresses this by redistributing flow through strategic edge swapping between surplus and saturated source-to-leaf paths.

\begin{lem}\longthmtitle{Rewire Identifies and Resolves Infeasibility}
    Given any capacity-infeasible solution from \texttt{FORWARD}, \textsf{Rewire} can identify blocking edges and redistribute flow to achieve feasibility when a feasible solution~exists.
    \label{lem::rewire_feasibility}
\end{lem}
\begin{pf}
    Given an unsupplied node $v$ (where $\vect{p}^\ell_v < 0$), by the flow balance enforced by \textsf{Islander}, there must exist a source-to-leaf path $\mathcal{R}_{j}^{i'}$ from polytree $\mathcal{T}^\ell_{i'}$ that can supply the remaining demand through a non-capacity-constrained route:
    {\small $$\exists \mathcal{R}_{j}^{i'} \in \mathcal{C}(\mathcal{R},v)\setminus\mathcal{C}(\bar{\mathcal{R}}) \text{ with leaf node } u \text{ such that } (u,v)\in\mathcal{G}^\ell$$}
    By flow balance, there also exists a polytree $\mathcal{T}^\ell_i$ with surplus flow, meaning $\mathcal{T}^\ell_i\cap\mathcal{C}(\mathcal{R}^+)\neq\emptyset$. This provides a source-to-leaf path $\mathcal{R}_{j}^{i}\in\mathcal{C}(\mathcal{R}^+)$ with leaf node $v'$ that has excess supply. \textsf{Rewire} performs the following swap: (1) connects the surplus leaf $v'$ to polytree $\mathcal{T}^\ell_{i'}$, (2) redirects the non-constrained path to supply node $v$, and (3) removes any blocking edge $\hat{e}_{u',v'}$ that would create a cycle. This maintains radiality while achieving feasibility.\boxend
\end{pf}
\begin{thm}\longthmtitle{\texttt{FORWARD} with \textsf{Rewire} Generates Feasible Solutions}
    The complete \texttt{FORWARD} algorithm including \textsf{Rewire} is guaranteed to construct feasible radial configurations for optimization problem \eqref{eqn::problem1} when feasible solutions exist.
    \label{thm::feasibility_complete}
\end{thm}
\begin{pf}
    The proof follows from Theorem~\ref{thm::feasibility_uncapacitated} and Lemma~\ref{lem::rewire_feasibility}. When capacity constraints are not binding, Theorem~\ref{thm::feasibility_uncapacitated} guarantees feasibility. When capacity constraints create infeasibilities, Lemma~\ref{lem::rewire_feasibility} demonstrates that \textsf{Rewire} can identify and resolve these infeasibilities through strategic flow redistribution while maintaining the radial structure. Since \textsf{Rewire} operates under the assumption that a feasible solution exists (Assumption~\ref{assump::feas}), the complete algorithm is guaranteed to find such a~solution.\boxend
\end{pf}

\section{Complexity} 
\label{subsec::complexity}
\vspace{-0.1in}

The computational complexity of \texttt{FORWARD} is determined by the sequential execution of its sub-processes. For each sub-process, we analyze their complexity as follows:

\begin{itemize}
\item \textsf{Pre-Processor}: $\mathcal{O}(n + m)$ as it processes all pendant nodes and their incident edges, potentially iterating through the entire graph structure.

\item \textsf{Islander}: $\mathcal{O}(n + m)$ for finding articulation points and partitioning the graph using standard graph algorithms.

\item \textsf{Net-Concad}: $\mathcal{O}(m)$ per call for condensing the graph into super-nodes using connected component algorithms.

\item \textsf{Sampler}: $\mathcal{O}(m \log m)$ per call due to priority queue operations over candidate edges. This function is called at most $n-1$ times to construct the spanning forest.

\item \textsf{Rewire}: $\mathcal{O}(n + m)$ in the worst case, where it may need to examine all source-to-leaf paths and perform edge swapping operations.
\end{itemize}

The overall complexity is dominated by the iterative sampling process: $\mathcal{O}((n + m) + (n + m) + (n-1) \cdot (m + m \log m) + (n + m))$, which simplifies to $\mathcal{O}(nm \log m)$.

However, this analysis represents the worst-case scenario. In practice, the complexity is significantly reduced due to several factors: \emph{Graph Sparsity:} Distribution networks typically exhibit sparse, small-world properties where $m = \mathcal{O}(n)$~\cite{BH-VS:20}. Under this assumption, the complexity reduces to $\mathcal{O}(n^2 \log n)$; \emph{Problem Decomposition:} After \textsf{Pre-Processor}, the remaining 2-core subgraph $\mathcal{G}_P$ contains significantly fewer than $n$ nodes. Subsequently, \textsf{Islander} partitions this reduced graph into even smaller subgraphs $\mathcal{G}^\ell$, each processed independently and potentially in parallel; \emph{Early Termination:} The \textsf{Sampler} often terminates before examining all $m$ edges due to the priority-based selection mechanism and the quasi-bipartite structure of the dual graph $\bar{\mathcal{G}}^\ell$.

Therefore, while the theoretical worst-case complexity is $\mathcal{O}(n^2 \log n)$ for sparse networks, the practical performance is substantially better due to the algorithm's structure-aware design and the inherent properties of distribution networks.


\section{Demonstration Example}
\label{sec::numerical}
\vspace{-0.1in}
To demonstrate the computational benefits of \texttt{FORWARD}, we evaluate its performance on optimal power distribution problems where generating feasible radial configurations is a computationally challenging and recurring task. We consider six different power distribution networks: three IEEE standard test systems (IEEE 13 with 2 sources, IEEE 18 with 2 sources, IEEE 33 with 3 sources) and three small-world networks (WS 120 with 10 sources, WS 240 with 10 sources, and WS 400 with 20 sources).

\subsection{Problem Formulation for Power Networks}
\vspace{-0.1in}
The IEEE graphs have been preprocessed to single-phase DC network structures. As established in the literature~\cite{russell}, single-phase direct current (DC) simplifications provide effective approximations to three-phase alternating current (AC) flow problems, allowing the optimal power flow problem to be cast in the form of Problem~\eqref{eqn::problem1}. The small-world networks are generated using the Watts-Strogatz mechanism~\cite{Watts1998} to create realistic network topologies with small-world properties.

In power distribution networks, the parameters of Problem~\eqref{eqn::problem1} have specific physical interpretations:
\begin{itemize}
\item $d_i$ represents power demand at load nodes (sinks)
\item $g_i$ represents power generation at generator nodes (sources)  
\item $x_{i,j}$ represents power flow through transmission line $(i,j)$.
\end{itemize}
The power loss across edge $e_{i,j}$ is computed as $R_{i,j} \left(\frac{x_{i,j}}{v_j}\right)^2$, where $v_j$
 is voltage at node $j$ and is a fix value. The voltage variation across the distribution network often is neglected. Using per unit ( p.u.) representation we assume that $v_j=1$ p.u. for all nodes $i\in\mathcal{V}_D$~\cite{QP-HL:13}. Therefore, in problem~\eqref{eqn::problem1} we have $C_{i,j}=R_{i,j}.$ 

For this problem we design the sampling weights based on how power flows across a distribution line. Following Ohm's law, the power flow $x_{i,j}$ through edge $(i \to j)$ depends on the total power demand $d_{i \to j}$ of all loads supplied through that edge and the voltage $v_j$ at node $j$. The power loss cost through edge $(i \to j)$ is:
$$f_{(i \to j)} = R_{i,j} \cdot \left(\frac{d_{i \to j}}{v_j}\right)^2$$

Since all edges within a source-to-leaf path $\mathcal{R}_k^i$ carry the same flow, we define a rolling-demand estimator:
$$\hat{f}(\mathcal{R}_k^i) = \sum\nolimits_{\kappa \in \mathcal{R}_k^i} R_{\kappa,\kappa+1} \cdot d_\kappa$$

The sampling weight for the \textsf{Sampler} function is then defined as:
\begin{equation}
    \label{eqn::distribution_prob2}
    w_{i,j} = \frac{\vect{p}_i}{R_{i,j} \cdot d_j^2 + \hat{f}(\mathcal{R}_k^i)}
\end{equation}
where $\vect{p}_i$ is the current nodal value representing remaining power at node $i$, which decreases at each iteration as power demand is allocated. The weights $w_{i,j}$ are dynamically updated and normalized at each iteration across all candidate edges.

\subsection{Experimental Setup, Results and Analysis}
\vspace{-0.1in}
We implemented our numerical experiments using the PowerDistributionModel (PMD) framework~\cite{PMD} from Los Alamos National Laboratory, coded in Julia. For comparison, we use Knitro from Artelys to solve the MINLP formulation of Problem~\eqref{eqn::problem1}. All experiments were conducted on a MacBook Air with an M3 chip and 24 GB of RAM. The source codes for this simulation study is available at~\cite{forward_github}.

The simulation results are presented in Table~\ref{tab:sample_table} and Figures~\ref{fig:values} and~\ref{fig:cpu}. Power loss values are reported in kilowatts (kW). In the time column, $^\star$ indicates processes manually terminated after excessive computation time, and ``TL" (Time Limit) indicates cases where no solution was found within 3 hours. The results demonstrate several key advantages of \texttt{FORWARD}:

\begin{figure*}[htbp]
    \centering
    \begin{minipage}{0.45\textwidth}
        \centering
        \includegraphics[width=0.9\textwidth]{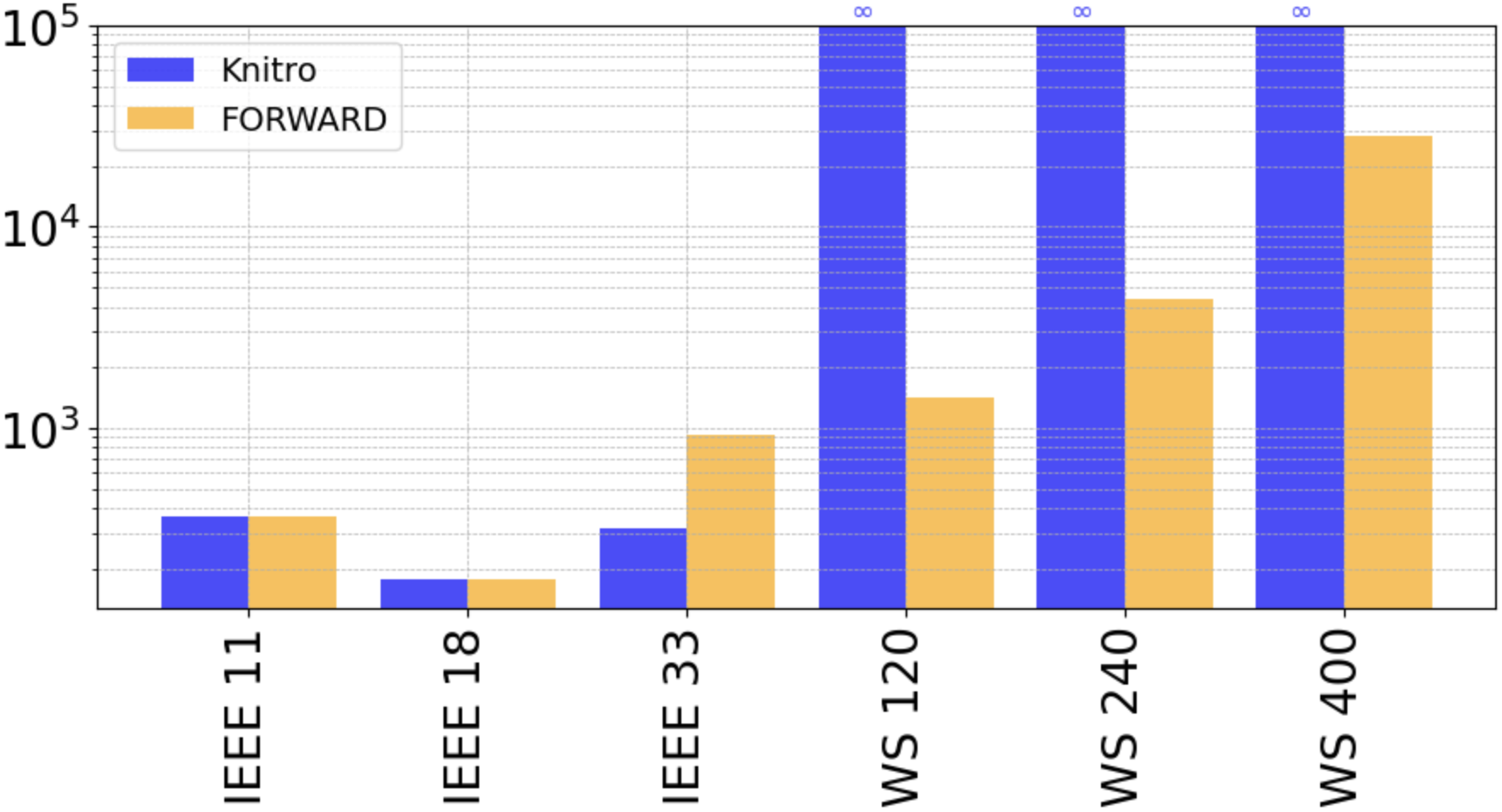}
        \caption{{\footnotesize Power loss comparison between Knitro and \texttt{FORWARD}.}}
        \label{fig:values}
    \end{minipage}\quad
    \begin{minipage}{0.45\textwidth}
        \centering
        \includegraphics[width=1\linewidth]{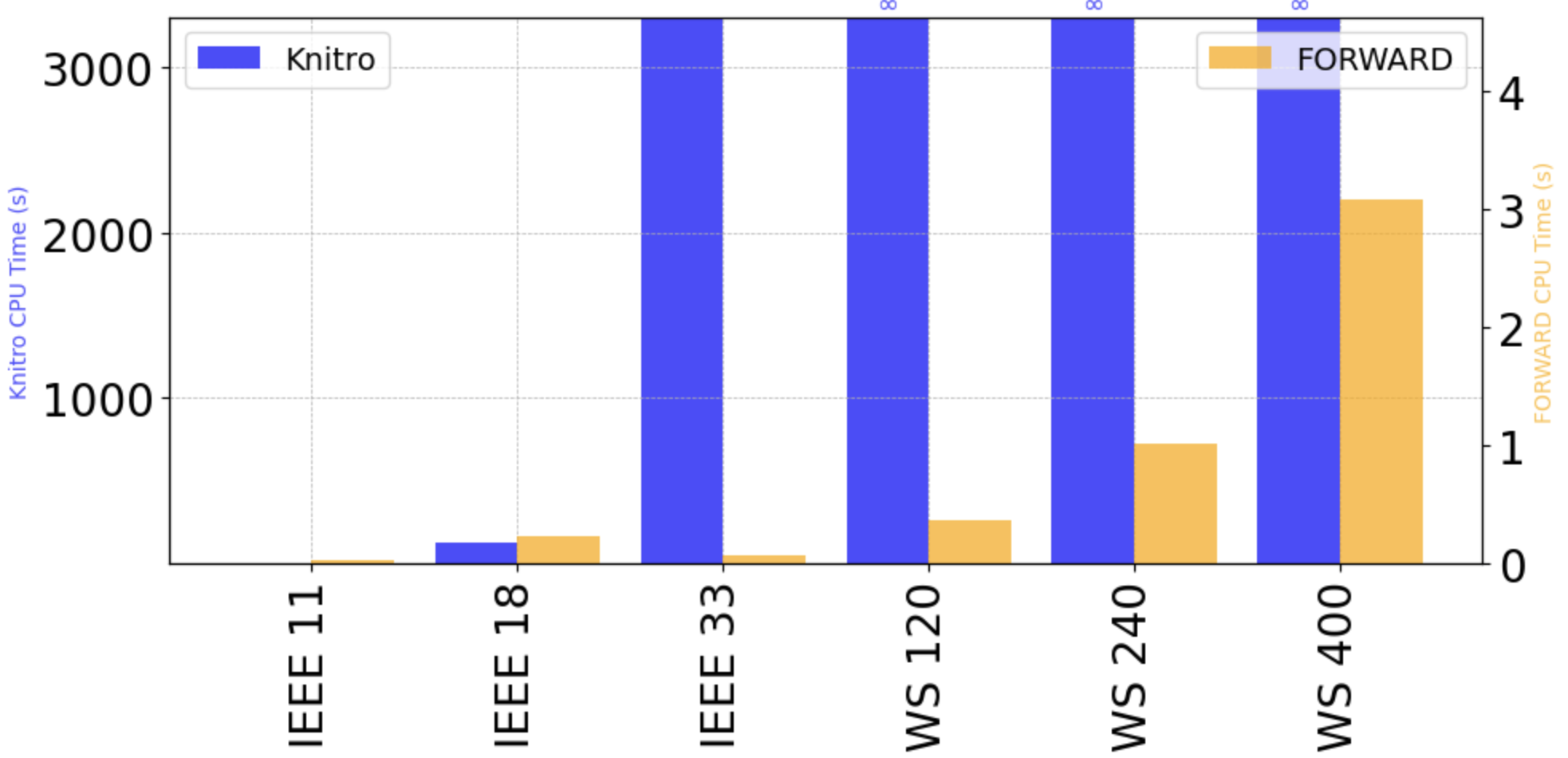}
        \caption{{\footnotesize CPU time comparison between Knitro and \texttt{FORWARD}.}}
        \label{fig:cpu}
    \end{minipage}
\end{figure*}

\begin{table*}[t]
\centering
\caption{{Performance comparison between Knitro and \texttt{FORWARD}.}}
\label{tab:sample_table}
{\scriptsize
\begin{tabular}{|c||c|c||c|c|}
\hline
~ & \multicolumn{2}{c||}{Knitro} & \multicolumn{2}{c|}{\texttt{FORWARD}} \\
\hline
Network & Power Loss (kW) & CPU Time (s) & Power Loss (kW) & CPU Time (s) \\ 
\hline
IEEE 13 & 360.183 & 4.189 & 360.183 & 0.033 \\ 
\hline
IEEE 18 & 175.821 & 123.416 & 175.821 & 0.229 \\ 
\hline
IEEE 33 & 318.568 & 3,345.67$^\star$ & 919.783 & 0.066 \\ 
\hline
WS 120 & TL & TL & 1,428.72 & 0.361 \\ 
\hline
WS 240 & TL & TL & 4,393.17 & 1.016 \\ 
\hline
WS 400 & TL & TL & 28,345.7 & 3.090 \\ 
\hline
\end{tabular}
}
\end{table*}

\emph{Computational Efficiency:} Knitro was only able to return solutions for small-sized networks, with CPU times significantly exceeding those of \texttt{FORWARD}. Commercial solvers like Knitro rely on heuristics for warm-start initialization within the feasible domain. As network size increases, these heuristics struggle to find proper initialization points within polynomial time due to the exponential growth of the combinatorial space, often resulting in prohibitive computation times or complete failure to find initial feasible points.

\emph{Scalability:} In contrast, \texttt{FORWARD} successfully constructed feasible radial configurations for all test cases with remarkably low computation times, even for networks with 400 nodes. The algorithm's polynomial-time complexity enables it to handle large-scale networks that are intractable for exact MINLP solvers.

\emph{Solution Quality:} For cases where both methods found solutions (IEEE 13 and IEEE 18), \texttt{FORWARD} achieved identical power loss values, demonstrating optimal performance on these instances. For the IEEE 33 network, while Knitro was terminated before completion, \texttt{FORWARD} provided a feasible solution with reasonable power loss in a fraction of the time. The promising solution quality can be attributed to the physically-aware weight design in Equation~\eqref{eqn::distribution_prob2} and the systematic sampling mechanism that considers both electrical properties and network topology.

\emph{Practical Impact:} These results highlight \texttt{FORWARD}'s practical value for real-world power distribution planning, where computational tractability is essential for operational decision-making. The algorithm's ability to handle networks of realistic size (400+ nodes) within seconds makes it suitable for online optimization and real-time grid management applications.


\section{Conclusions}
\label{sec:conclusions}
\vspace{-0.1in}
This paper presented \texttt{FORWARD}, a polynomial-time algorithm for constructing feasible radial configurations in multi-source distribution networks. We addressed the optimal radial reconfiguration problem~\eqref{eqn::problem1} and proved that feasibility determination was weakly NP-complete, making exact methods computationally intractable for large networks. Our approach consisted of three main contributions. First, we developed a graph decomposition strategy (Theorem~\ref{thm::partition}) that enabled parallel processing while preserving feasibility and optimality. Second, we designed a five-function algorithmic framework where \textsf{Net-Concad} addressed greedy shortsightedness through dual graph condensation, and \textsf{Rewire} handled capacity constraints via strategic edge swapping. We provided rigorous theoretical analysis proving feasibility guarantees for both uncapacitated and capacitated cases. Third, numerical evaluation on six networks (up to 400 nodes) demonstrated \texttt{FORWARD}'s superior computational efficiency compared to commercial MINLP solvers. While traditional methods required hours or failed entirely, \texttt{FORWARD} achieved solutions in seconds with optimal or near-optimal quality. The algorithm's success stemmed from leveraging radial network structure and incorporating physically-motivated weight functions rather than relying on generic branch-and-bound approaches. The polynomial-time complexity made \texttt{FORWARD} suitable for real-time distribution network management applications. Despite being designed for quadratic cost models, the algorithm successfully handled complex physics constraints in power system reconfiguration, establishing its effectiveness as an initialization strategy for iterative optimization solvers.

Future work included formally characterizing the optimality gap, extending to nonlinear constraint models, and investigating applicability to broader network optimization problems. The promising results suggested that our algorithmic principles could impact telecommunications, supply chain, and transportation planning beyond radial distribution networks.

\bibliographystyle{plain} 
\bibliography{references}  

\appendix

\section{Proofs}\label{Appedix::Proofs}
\vspace{-0.1in}
\begin{pf}[Proof of Theorem~\ref{thm::partition}]
Let $\mathcal{G}^\ell = (\mathcal{V}^\ell, \mathcal{E}^\ell)$ denote the $\ell$-th component. For any separation node $v$, let $\mathcal{C}(v) = \{\ell | v \in \mathcal{V}^\ell\}$ denote the set of components containing $v$.

\textbf{($\Rightarrow$)} Suppose there exists a feasible solution $(\mathcal{S}, \vect{x})$ for problem~\eqref{eqn::problem1} on $\mathcal{G}_D$. 

For each component $\mathcal{G}^\ell$ with separation nodes $\mathcal{V}_{sep}^\ell = \{v_1, v_2, \ldots, v_k\} \subset \mathcal{V}^\ell$, we construct the adjusted nodal values as follows: For each separation node $v_i \in \mathcal{V}_{sep}^\ell$, Kirchhoff's law in the original graph gives:
$$\sum_{u \in \mathcal{V}^\ell \setminus \{v_i\}} (x_{u,v_i} - x_{v_i,u}) + \sum_{\ell' \neq \ell} \sum_{u \in \mathcal{V}^{\ell'}} (x_{u,v_i} - x_{v_i,u}) = g_{v_i} - d_{v_i}.$$
Define:
\begin{align*}
f_{v_i}^{\ell}
  &= \sum_{u \in \mathcal{V}^\ell \setminus \{v_i\}} (x_{u,v_i} - x_{v_i,u})
  && \text{\parbox[t]{.45\linewidth}{(net flow into $v_i$ from within $\mathcal{G}^\ell$)}}\\
f_{v_i}^{\text{ext}}
  &= \sum_{\ell' \neq \ell} \sum_{u \in \mathcal{V}^{\ell'}} (x_{u,v_i} - x_{v_i,u})
  && \text{\parbox[t]{.45\linewidth}{(net external flow into $v_i$)}}.
\end{align*}
Then $f_{v_i}^{\ell} + f_{v_i}^{ext} = g_{v_i} - d_{v_i}$ for each $i \in \{1,2,\ldots,k\}$. Define the adjusted nodal values:
\begin{align*}
p_{v_i}^{\ell} &= g_{v_i} - d_{v_i} - f_{v_i}^{ext} = f_{v_i}^{\ell} \quad \text{for separation nodes}\\
p_u^{\ell} &= g_u - d_u \quad \text{for internal nodes } u \in \mathcal{V}^\ell \setminus \mathcal{V}_{sep}^\ell.
\end{align*}
To verify that component $\mathcal{G}^\ell$ is balanced, we compute:
$$\sum\nolimits_{u \in \mathcal{V}^\ell} p_u^{\ell} = \sum\nolimits_{u \in \mathcal{V}^\ell \setminus \mathcal{V}_{sep}^\ell} (g_u - d_u) + \sum\nolimits_{i=1}^k f_{v_i}^{\ell}.$$
By Kirchhoff's law applied to the subgraph $\mathcal{G}^\ell$, since internal nodes can only exchange flow with separation nodes (no external connections):
$$\sum\nolimits_{u \in \mathcal{V}^\ell \setminus \mathcal{V}_{sep}^\ell} (g_u - d_u) + \sum\nolimits_{i=1}^k f_{v_i}^{\ell} = 0.$$
Therefore, $\sum_{u \in \mathcal{V}^\ell} p_u^{\ell} = 0$, confirming that $\mathcal{G}^\ell$ is balanced. When we restrict the original flow $\vect{x}$ to edges within $\mathcal{G}^\ell$, obtaining $\vect{x}^\ell$, the flow conservation is automatically satisfied:
\begin{itemize}
\item For internal node $u \in \mathcal{V}^\ell \setminus \mathcal{V}_{sep}^\ell$: 
$$\sum\nolimits_{w|(u,w) \in \mathcal{E}^\ell} (x_{w,u}^\ell - x_{u,w}^\ell) = g_u - d_u = p_u^\ell$$
\item For separation node $v_i \in \mathcal{V}_{sep}^\ell$: 
$$\sum\nolimits_{w|(v_i,w) \in \mathcal{E}^\ell} (x_{w,v_i}^\ell - x_{v_i,w}^\ell) = f_{v_i}^{\ell} = p_{v_i}^\ell$$
\end{itemize}

Hence, $A(\mathcal{S}^\ell)\vect{x}^\ell = \vect{p}^\ell$ is satisfied for each component $\mathcal{G}^\ell$.
\textbf{($\Leftarrow$)} Suppose there exist feasible solutions $(\mathcal{S}^\ell, \vect{x}^\ell)$ for each balanced component $\mathcal{G}^\ell$ satisfying $A(\mathcal{S}^\ell)\vect{x}^\ell = \vect{p}^\ell$.

Construct the global solution as $\mathcal{S} = \bigcup_{\ell=0}^L \mathcal{S}^\ell$ and $\vect{x}$ as the concatenation of all $\vect{x}^\ell$.

The capacity constraints are satisfied since they hold for each component individually and the edge sets are disjoint. The union of disjoint polyforests rooted at source nodes remains a polyforest, preserving the radial property.

For flow conservation at any separation node $v$ shared by multiple components, we have:
$$\sum\nolimits_{\ell \in \mathcal{C}(v)} \left(\sum\nolimits_{u \in \mathcal{V}^\ell \setminus \{v\}} (x_{u,v}^\ell - x_{v,u}^\ell)\right) = \sum\nolimits_{\ell \in \mathcal{C}(v)} p_v^{\ell}.$$
By construction of the adjusted values:
$$\sum\nolimits_{\ell \in \mathcal{C}(v)} p_v^{\ell} = \sum\nolimits_{\ell \in \mathcal{C}(v)} (g_v - d_v - f_v^{ext,\ell}) = g_v - d_v$$
where the external flow terms cancel out since $\sum_{\ell \in \mathcal{C}(v)} f_v^{ext,\ell}$ $= 0$ (what flows out of one component flows into another at separation node $v$). For internal nodes (non-separation), flow conservation is trivially satisfied since their nodal values remain unchanged. Therefore, Kirchhoff's law is satisfied at every node in the reconstructed graph, completing the proof.
\boxend\end{pf}
\begin{pf}[Proof of Corollary~\ref{cor::optimality_partition}]
The proof follows from the separability of the quadratic cost function across disjoint edge sets and the equivalence established in Theorem~\ref{thm::partition}. Since the components $\mathcal{G}^\ell$ have disjoint edge sets $\mathcal{E}^\ell$, the total cost of any feasible solution $(\mathcal{S}, \vect{x})$ on $\mathcal{G}_D$ can be decomposed as:
$$\sum\nolimits_{(i,j)\in\mathcal{S}} C_{i,j} \cdot x_{i,j}^2 = \sum\nolimits_{\ell=0}^L \sum\nolimits_{(i,j)\in\mathcal{S}^\ell} C_{i,j} \cdot (x_{i,j}^\ell)^2,$$
where $\mathcal{S}^\ell = \mathcal{S} \cap \mathcal{E}^\ell$ and $x_{i,j}^\ell$ is the flow on edge $(i,j)$ within component $\mathcal{G}^\ell$. By Theorem~\ref{thm::partition}, there is a bijective correspondence between feasible solutions on $\mathcal{G}_D$ and collections of feasible solutions on the balanced components $\{\mathcal{G}^\ell\}_{\ell=0}^L$. Therefore:
$$\min_{(\mathcal{S}, \vect{x}) \text{ feasible on } \mathcal{G}_D} \sum\nolimits_{(i,j)\in\mathcal{S}} C_{i,j} \cdot x_{i,j}^2 = $$$$\sum\nolimits_{\ell=0}^L \min_{(\mathcal{S}^\ell, \vect{x}^\ell) \text{ feasible on } \mathcal{G}^\ell} \sum\nolimits_{(i,j)\in\mathcal{S}^\ell} C_{i,j} \cdot (x_{i,j}^\ell)^2$$

Since $(\mathcal{S}^{\ell,*}, \vect{x}^{\ell,*})$ is optimal for component $\mathcal{G}^\ell$, the combined solution achieves this minimum cost. By Theorem~\ref{thm::partition}, this combined solution is feasible for the original problem, and by the cost separability, it is optimal.
\boxend\end{pf}

\section{Sub-routines}
\label{Appedix::subroutines}

In this appendix we introduce with greater detail internal sub-routines in \texttt{FORWARD} mentioned in Section~\ref{sec::algorithm}.

{\small
\begin{algorithm} 
{\footnotesize
\caption{\textsf{Tree-Update}} 
\begin{algorithmic}[1]
\Require Existing polytree set $\mathbb{T}^\ell$ and newly sampled edge $i\to j$.
\For {$\mathcal{T}\in\mathbb{T}^\ell$}
\If { $i\in\mathcal{V}(\mathcal{T})$}
\State $\mathbb{T}^\ell\gets \mathbb{T}^\ell\backslash \mathcal{T}$
\State $\mathcal{T}\gets \mathcal{T}\cup(\{i,j\},i\to j)$ \Comment{New edge added to polytree}
\For{$\mathcal{T}'\in\mathbb{T}^\ell\backslash\mathcal{T}$}  \Comment{Merging polytrees if needed}
\If { $j\in\mathcal{V}(\mathcal{T}')$} 
\State $\mathbb{T}^\ell\gets \mathbb{T}^\ell\backslash \mathcal{T}'$
\State $\mathcal{T}\gets \mathcal{T}\cup\mathcal{T}'$
\EndIf
\EndFor
\State \textbf{Break}
\EndIf 
\EndFor
\State \Return $\mathbb{T}^\ell$
\end{algorithmic}
\label{alg::tree_update}
}
\end{algorithm}
}

{\small
\begin{algorithm} 
{\footnotesize
\caption{\textsf{Conncomp}~\cite{conncomp}} 
\begin{algorithmic}[1]
\Require Graph $\mathcal{G} = (\mathcal{V}, \mathcal{E})$
\Ensure Set of connected components $\mathbb{C}$

\State $\mathbb{C} \gets \emptyset$
\State $\text{visited}[v] \gets \textbf{false}$ for all $v \in \mathcal{V}$
\For{each $v \in \mathcal{V}$}
    \If{not $\text{visited}[v]$}
        \State $\mathcal{C} \gets \emptyset$
        \State \Call{DFS}{$v$, $\mathcal{G}$, visited, $\mathcal{C}$}
        \State $\mathbb{C} \gets \mathbb{C} \cup \{\mathcal{C}\}$
    \EndIf
\EndFor
\State \Return $\mathbb{C}$
\vspace{1mm}
\Function{DFS}{$v$, $\mathcal{G}$, visited, $\mathcal{C}$}
    \State $\text{visited}[v] \gets \textbf{true}$
    \State $\mathcal{C} \gets \mathcal{C} \cup \{v\}$
    \For{each $u \in \mathcal{N}(v)$}
        \If{not $\text{visited}[u]$}
            \State \Call{DFS}{$u$, $\mathcal{G}$, visited, $\mathcal{C}$}
        \EndIf
    \EndFor
\EndFunction
\end{algorithmic}
\label{alg::connected_components}
}
\end{algorithm}
}

{\small
\begin{algorithm} 
{\footnotesize
\caption{\textsf{Retrieve}} 
\begin{algorithmic}[1]
\Require Existing queue $q = \{x_0,x_1,\cdots,x_n\}$.
\State $q \gets q\setminus\{x_0\}$
\State \Return $x_0$
\end{algorithmic}
\label{alg::retrieve}
}
\end{algorithm}
}

\end{document}